\documentclass[12pt,reqno]{amsart}%
\usepackage{amsmath}%
\usepackage{mathtools}
\usepackage{amssymb}%
\usepackage{amscd}%
\usepackage{color}%
\usepackage{enumitem}%
\usepackage{graphicx}%
\usepackage{mathrsfs}%
\usepackage[all]{xy}%
\usepackage{stmaryrd}
\usepackage{newtxtext,newtxmath}%
\usepackage[all, warning]{onlyamsmath}
\usepackage{amsrefs}
\usepackage{bbm}
\usepackage{tikz}
\usepackage{longtable}
\usepackage{array}
\usepackage[capitalise]{cleveref}
\def\wt#1{\widetilde{#1}}
\usepackage{aliascnt}
\newcommand{\newsharedtheorem}[3]{%
  \newaliascnt{#1}{#2}%
  \newtheorem{#1}[#1]{#3}%
  \aliascntresetthe{#1}%
  \crefname{#1}{#3}{#3s}%
}
\theoremstyle{plain}
    
    \newtheorem{theorem}{Theorem}[section]
    \newsharedtheorem{proposition}{theorem}{Proposition}
    \newsharedtheorem{lemma}{theorem}{Lemma}
    \newsharedtheorem{corollary}{theorem}{Corollary}
    \newsharedtheorem{corint}{thmint}{Corollary}
    \newsharedtheorem{problem}{theorem}{Problem}
     
\theoremstyle{definition}
    \newsharedtheorem{definition}{theorem}{Definition}
    \newsharedtheorem{notation}{theorem}{Notation}
    \newsharedtheorem{example}{theorem}{Example}
    \newsharedtheorem{remark}{theorem}{Remark}
\def\Alphabet{A,B,C,D,E,F,G,H,I,J,K,L,M,N,O,P,Q,R,S,T,U,V,W,X,Y,Z}
\def\grabet{a,b,c,d,e,f,g,h,i,j,k,l,m,n,o,p,q,r,s,t,u,v,w,x,y,z}
\def\endpiece{xxx}
\def\makeAlphabet[#1]{\expandafter\makeA#1,xxx,}
\def\makealphabet[#1]{\expandafter\makea#1,xxx,}
\def\makeA#1,{\def\temp{#1}\ifx\temp\endpiece\else%
\mkbb{#1}\mkfrak{#1}\mkbf{#1}\mkcal{#1}\mkscr{#1}\mkbs{#1}\expandafter\makeA\fi}%
\def\makea#1,{\def\temp{#1}\ifx\temp\endpiece\else\mkfrak{#1}\mkbf{#1}\mkbs{#1}\expandafter\makea\fi}%
\def\mkbb#1{\expandafter\def\csname bb#1\endcsname{\mathbb{#1}}}
\def\mkfrak#1{\expandafter\def\csname fr#1\endcsname{\mathfrak{#1}}}
\def\mkbf#1{\expandafter\def\csname b#1\endcsname{\mathbf{#1}}}
\def\mkcal#1{\expandafter\def\csname c#1\endcsname{\mathcal{#1}}}
\def\mkscr#1{\expandafter\def\csname s#1\endcsname{\mathscr{#1}}}
\def\mkbs#1{\expandafter\def\csname bs#1\endcsname{{\boldsymbol{#1}}}}
\def\makeop[#1]{\xmakeop#1,xxx,}
\def\mkop#1{\expandafter\def\csname #1\endcsname{{\mathrm{#1}}}} %
\def\xmakeop#1,{\def\temp{#1}\ifx\temp\endpiece\else\mkop{#1}\expandafter\xmakeop\fi}%
\def\makeup[#1]{\xmakeup#1,xxx,}
\def\mkup#1{\expandafter\def\csname #1\endcsname{{\mathrm{#1}\,}}} %
\def\xmakeup#1,{\def\temp{#1}\ifx\temp\endpiece\else\mkup{#1}\expandafter\xmakeup\fi}%
\makeAlphabet[\Alphabet]
\makealphabet[\grabet]
\makeop[Map,Hom,id,Fun,cont,Avg,Var,Ob,unif,loc,pr,End,Binom,ev]
\makeup[Im,Ker,Coker]
\makeop[Cov,Sym]
\def\Stoch{\mathbf{Stoch}}

\def\ve{\varepsilon}
\def\bspi{\boldsymbol{\pi}}

\def\R{\bbR}

\def\pair#1{\langle#1\rangle}
\def\abs#1{|#1|}

\def\dabs#1{\Vert#1\Vert}

\def\bsone{\boldsymbol{1}}

\def\one{\boldsymbol{1}}

\begin{document}

\setcounter{tocdepth}{1}
\newpage
\title{Bayesian Filters via the Stochastic Category}
\author[K. Bannai]{Kenichi Bannai}\email{bannai@math.keio.ac.jp}
\author[M. Kurisaki]{Masahiro Kurisaki}
\thanks{This work is supported in part by RIKEN AIP, as well as JST CREST Grant Number JPMJCR2112 (The Bayes-Duality Project FY2021--2026, led by Emtiyaz Khan).
The second author is supported by the RIKEN Special Postdoctoral Researcher Program, FY2025--2027.}
\address[Bannai]{Department of Mathematics, Faculty of Science and Technology, Keio University, 3-14-1 Hiyoshi, Kouhoku-ku, Yokohama 223-8522, Japan.}
\address[Bannai, Kurisaki]{Mathematical Science Team, RIKEN Center for Advanced Intelligence Project (AIP), 1-4-1 Nihonbashi, Chuo-ku, Tokyo 103-0027, Japan.}

\date{November 1, 2026}
\begin{abstract}
	We extend the categorical treatment of Bayesian filters to models with feedback,
	whose transitions and observations may depend on the past observations.
	We work in the \emph{stochastic category}, whose objects are measurable spaces and whose
	morphisms are \emph{Markov kernels},
	which we refer to as \emph{stochastic morphisms}.
	We systematically parameterize stochastic morphisms by a measurable space,
	and develop the calculus of conditioning in a notation intuitively 
	compatible with classical probability theory,
	requiring neither string diagrams nor densities.
	A stochastic morphism into a product space is regarded as a parameterized family of joint laws,
	and its conditional law is again a stochastic morphism,
	so that the rules of conditioning become equalities of stochastic morphisms.
	We define the Bayesian inverse with respect to a parameterized prior and prove Bayesian coherence.
	A model with feedback is then a model whose parameter space grows with the observations,
	and we derive the Bayesian filter, prediction and smoothing for such models,
	with standard Borel hidden states and countably generated observation spaces.
	The affine Gaussian stochastic morphisms parameterized by a measurable space
	form a family strictly larger than the Markov category \textbf{Gauss},
	closed under composition and Bayesian inversion,
	in which we obtain the Kalman filter and the Rauch--Tung--Striebel smoother with feedback.
\end{abstract}

\maketitle

%
%
%
\section{Introduction}
%
%
%

Stochastic morphisms, also known as Markov kernels, between measurable spaces
were introduced by Lawvere \cite{La62} to express stochastic phenomena
such as noise and uncertainty in categorical terms,
and independently by \v{C}encov \cite{C65,C82} from the viewpoint of statistical decision theory.
The resulting category, which we call the \emph{stochastic category}, 
was interpreted by Giry \cite{G82} as the Kleisli category of the probability monad,
and its structure was axiomatized by Fritz \cite{F20} as Markov categories (see also \cite{FGPR23}),
in which marginalization, conditioning and independence are formulated through compositions and monoidal products of morphisms.
In this framework, arguments are carried out with string diagrams,
and classical algorithms such as the Bayesian filter \cite{F25,V23} and the smoother \cite{F25} can be derived uniformly.

In classical probability theory, Bayesian filters and their continuous time analogues have been studied extensively
\cite{Kal60,CMR05,BC09,SS23,LS01}, and filters are naturally used with feedback, i.e.\
the transitions and observations of the model may depend on the observations made so far,
for instance through a control input or through coefficients chosen based on the data \cite[\S4.2]{CMR05}.
Such models are not hidden Markov models in the sense of \cite{F25}, since the future is not conditionally independent of the past given the current hidden state.
They are instead characterized by their conditional laws given the whole past (\cref{prop: characterization}).
The purpose of this article is to extend the categorical interpretation of Bayesian filters to this setting,
by systematically parameterizing stochastic morphisms,
using notation intuitively compatible with classical probability.
In categorical terms, fixed parameters are treated in \cite{F25} through the parametric Markov categories
of \cite[\S2.2]{FGPR23}, and inputs appear in \cite[Definition~3.1]{V23}.
Parameterized stochastic morphisms
also appear in \cite{For21} as transitional random variables.
In this article, the parameter space grows with the observations.
For simplicity, we work in the stochastic category, the motivating example of a Markov category.

We take the joint law as the primary object.
A joint law is a stochastic morphism $\pi\colon S\rightarrow X\times Y$,
that is, a family of probability measures $\pi(s)$ on $X\times Y$ parameterized by a measurable space $S$.
We write
\[
	\pi(Y\mid X)\colon S\times X\rightarrow Y
\]
for the conditional law of $Y$ given $X$ for $\pi$ (see \cref{def: CL}), a family of regular conditional probabilities
for the measures $\pi(s)$, characterized by
\[
	\pi=(\pr_X,\pi(Y\mid X))\circ(\id_S,\pi_X).
\]
The notation makes explicit the joint law and the parameter space,
and the factors $X$ and $Y$ play the role traditionally played by random variables,
conditioning on $X$ being conditioning on the projection $\pr_X$ (see \cref{rem: prX}).
The essential point is that $\pi(Y\mid X)$ is again a stochastic morphism, with parameter space $S\times X$,
so that further conditional laws can be taken and the rules of conditioning become equalities of stochastic morphisms.
For example, the tower property (\cref{thm: tower}) for $\pi\colon S\rightarrow X\times Y\times Z$ reads
\[
	\pi(Y\times Z\mid X)(Z\mid Y)=\pi(Z\mid X\times Y),
\]
that is, conditioning on $X$ and then on $Y$ is the same as conditioning on $X\times Y$ at once,
corresponding to the identity $p(z\mid x,y)=p(y,z\mid x)/p(y\mid x)$ for densities.
No reference measure or density is required, and measure theory is confined to the existence and uniqueness
of the conditional law.

To see the calculus at work, consider the simplified case of a hidden Markov model without feedback
when $S$ is a single point, with initial law $\mu$ on $X_1$, transitions $a_k\colon X_k\rightarrow X_{k+1}$ and observations $h_k\colon X_k\rightarrow Y_k$ for integers $k\geq 1$.
Let $b_k\coloneqq(h_{k+1},\id_{X_{k+1}})\circ a_k\colon X_k\rightarrow Y_{k+1}\times X_{k+1}$,
which generates the next observation $Y_{k+1}$ and hidden state $X_{k+1}$ from $X_k$.
The joint law $\pi_k$ of the observations $Y^{[k]}\coloneqq Y_1\times\cdots\times Y_k$ and the hidden state $X_k$ 
is the probability measure on $Y^{[k]}\times X_k$ given by push-forwards
\[
	\pi_1\coloneqq(h_1,\id_{X_1})(\mu),\qquad \pi_{k+1}\coloneqq(\id_{Y^{[k]}}\times b_k)(\pi_k)
\]
of measures, and the posterior is the conditional law
\[
	\pi_k(X_k\mid Y^{[k]})\colon Y^{[k]}\rightarrow X_k
\]
giving a family of probability measures on $X_k$
parameterized by the observations $Y^{[k]}$ up to time $k$.
The following two rules of the calculus give the Bayesian filter.
Push-forward commutes with conditioning, so that
\[
	\pi_{k+1}(Y_{k+1}\times X_{k+1}\mid Y^{[k]})=b_k\circ\pi_k(X_k\mid Y^{[k]}),
\]
and the tower property then gives
\[
	\pi_{k+1}(X_{k+1}\mid Y^{[k+1]})=(b_k\circ\pi_k(X_k\mid Y^{[k]}))(X_{k+1}\mid Y_{k+1}).
\]
Since $b_k\circ\pi_k(X_k\mid Y^{[k]})=(h_{k+1},\id_{X_{k+1}})\circ\mu_{k+1}$
for the predicted prior $\mu_{k+1}\coloneqq a_k\circ\pi_k(X_k\mid Y^{[k]})\colon Y^{[k]}\rightarrow X_{k+1}$,
the right hand side is the Bayesian inverse $h_{k+1}^{\mu_{k+1}}$ of $h_{k+1}$ with respect to $\mu_{k+1}$,
a parameterized version of the Bayesian inverse of Cho and Jacobs \cite{CJ19}.
This gives the Bayesian filter
\[
	\pi_{k+1}(X_{k+1}\mid Y^{[k+1]})=h_{k+1}^{\mu_{k+1}},
\]
which states that the posterior, a priori involving the whole model up to time $k+1$,
is the inverse of a single observation $h_{k+1}$, the past entering only through the predicted prior $\mu_{k+1}$.

The parameter space enables the iteration.
The posterior $\pi_k(X_k\mid Y^{[k]})$ is a stochastic morphism with parameter space $Y^{[k]}$,
and serves as the prior for the next observation.
It also accommodates feedback.
A transition or observation depending on the past observations is a stochastic morphism with parameter space $Y^{[k]}$,
so that a model with feedback is simply a model whose parameter space grows with the observations.
Our main results, \cref{thm: BF,thm: PS}, establish the Bayesian filter for the joint law $\bspi$ of the model with feedback,
together with the corresponding recursions for prediction and smoothing,
every step being a Bayesian inversion with the past observations entering as the parameter.
In the case without feedback, the filter and the smoother are those of \cite{F25},
obtained in Markov categories with conditionals, where the past observations are carried as an input of the conditional.

The hidden states are standard Borel, the observation spaces only countably generated, and the parameter space arbitrary.
Since we work directly with conditional laws rather than conditional expectations,
no integrability is assumed, and the results are equalities of stochastic morphisms
holding for every value of the parameter $s\in S$.
In the affine Gaussian case with feedback, the filter and the smoother are classical
for conditionally Gaussian sequences \cite[Chapter~13]{LS01}, under integrability assumptions on the coefficients, and the filter is known to require only that the coefficients be almost surely finite \cite{CKvS89}.
In this article, the filter and the smoother follow from \cref{thm: BF,thm: PS} for arbitrary measurable coefficients, in a family of affine Gaussian stochastic morphisms strictly larger than the Markov category \textbf{Gauss} of \cite[\S6]{F20}.
Feedback is also natural in continuous time filtering, which will be treated in a subsequent article.

In classical terms, the conditional expectation of a bounded measurable function $f$ on $Y$ given $X$
is the pull-back
\[
	\bbE_\pi[f\mid X](s,x)=\int_Yf(y)\,\pi(Y\mid X)(s,x)(dy),\quad (s,x)\in S\times X,
\]
and the tower property of conditional expectations follows from that of conditional laws (\cref{prop: tower}).
In this sense, the calculus of conditional laws gives the calculus of conditional expectations.

The article is organized as follows.
In \S\ref{sec: CL} we first introduce conditional laws and augmentations.
We then prove fundamental results
such as the tower property,
and formulate conditional independence.
In \S\ref{sec: BI} we define the Bayesian inverse and prove Bayesian coherence:
updating observations successively gives the same posterior as updating them at once.
In \S\ref{sec: filter} we treat filtering, prediction and smoothing for hidden Markov models with feedback,
and in \S\ref{sec: Gaussian} the Kalman filter and the Rauch--Tung--Striebel smoother.
Appendix \ref{sec: appendix A} gives a self-contained proof of the existence and uniqueness of the conditional law,
which may be read independently,
Appendix \ref{sec: appendix B} relates our framework to classical expectations and conditional expectations,
and Appendix \ref{sec: appendix C} gives a dictionary between classical notation, the notation of Markov categories and our article.

%
%
%
\section{Conditional Laws and Augmentations}\label{sec: CL}
%
%
%

Throughout this article, $S,X,Y,Z,V,W$ denote measurable spaces.
In this section, we define the conditional law of a parameterized joint law
and prove fundamental properties such as the tower property, using augmentations of stochastic morphisms.

A measurable space $X$ is \emph{countably generated} if its $\sigma$-algebra is generated by
a countable family of measurable subsets, and \emph{standard Borel} if it is isomorphic as a measurable space
to a Borel subset of a Polish space with its Borel $\sigma$-algebra, or equivalently to a Borel subset of $[0,1]$ (see
\cite[Theorem~3.3.13]{Sri98}).
A standard Borel space is countably generated, and finite or countable products of 
standard Borel spaces are again standard Borel spaces.

A \emph{stochastic morphism} $p\colon X\rightarrow Y$, also known as a Markov kernel, probabilistic mapping
or transition probability from $X$ to $Y$,
 is a map sending any $x\in X$ to a
probability measure $p(x)$ on $Y$, such that for any measurable $B\subset Y$,
the map $x\mapsto p(x)(B)$ is measurable in $x$.
Any measurable map $f\colon X\rightarrow Y$ can be viewed
as a stochastic morphism by mapping $x$ to the Dirac measure $\delta_{f(x)}$ on $Y$.
We will call such a stochastic morphism a \emph{deterministic morphism}.
We denote by $\bsone\coloneqq\{*\}$ the one-point measurable space.
A probability measure $\mu$ on $X$ defines a stochastic morphism $\mu\colon\bsone\rightarrow X$
mapping $*$ to $\mu$,
and we will identify the two.

Given stochastic morphisms
$p\colon S\rightarrow X$ and $q\colon X\rightarrow Y$,
we define the composition by
\[
	(q\circ p)(s)(B)\coloneqq\int_Xq(x)(B)p(s)(dx)
\]
for any $s\in S$ and measurable set $B\subset Y$.
Then $q\circ p$ is again a stochastic morphism.
Measurable spaces and stochastic morphisms form a category, the \emph{stochastic category},
denoted $\Stoch$ in \cite{F20}, with the deterministic morphisms induced by the identities as identities.

We denote by $\pr_X\colon X\times Y\rightarrow X$ and
$\pr_Y\colon X\times Y\rightarrow Y$ the deterministic morphisms
induced by the projections.
More generally, we take the convention that on a product with several factors including $X$ and $Y$,
we denote by $\pr_{X\times Y}$ and $\pr_{Y\times X}$ the projections to $X\times Y$
and to $Y\times X$, the order of the indices indicating the order of the factors
in the range, and similarly for projections to more than two factors.
For stochastic morphisms $p\colon S\rightarrow X$ and $q\colon S\rightarrow Y$, the product
\[
	(p,q)\colon S\rightarrow X\times Y
\]
is the stochastic morphism $s\mapsto p(s)q(s)$,
where $p(s)q(s)$ is the product measure on $X\times Y$ given by $p(s)q(s)(A\times B)=p(s)(A)q(s)(B)$
for any $s\in S$ and measurable sets $A\subset X$ and $B\subset Y$.
For stochastic morphisms $a\colon X\rightarrow X'$ and $b\colon Y\rightarrow Y'$, 
we let $(a\times b)(x,y)\coloneqq a(x)b(y)$ for any $(x,y)\in X\times Y$.
We freely use the identities
\begin{equation}\label{eq: product}
	(a\times b)\circ(p,q)=(a\circ p,\,b\circ q),\qquad (p,q)\circ u=(p\circ u,\,q\circ u),
\end{equation}
the latter for a deterministic morphism $u$.
Furthermore, for the projection $\pr_X\colon X\times Y\rightarrow X$,
we have $\pr_X\circ(p,q)=p$.

We call a stochastic morphism $\pi\colon S\rightarrow X\times Y$ a \emph{joint law} on $X\times Y$ parameterized by $S$, or simply a joint law.
This gives for each $s\in S$ a probability measure $\pi(s)$ on $X\times Y$.
We define the \emph{marginal} of $\pi$ at $X$ by $\pi_X\coloneqq\pr_X\circ\pi$.
We define the joint laws and marginals similarly for any product of factors and subproduct.

\begin{definition}\label{def: CL}
	Consider a joint law $\pi\colon S\rightarrow X\times Y$.
	A \emph{conditional law} $\pi(Y\mid X)$ of $Y$ given $X$ for $\pi$ is a stochastic morphism
	\[
		\pi(Y\mid X)\colon S\times X\rightarrow Y
	\]
	that satisfies
	\begin{equation}\label{eq: CL}
		\pi=(\pr_X,\pi(Y\mid X))\circ(\id_S,\pi_X),
	\end{equation}
	where $\pi_X=\pr_X\circ\pi$ is the marginal.
\end{definition}

It is important to note that the $\pi$ in the notation $\pi(Y\mid X)$ is not a generic symbol to denote a
conditional law, but signifies that the conditional law is induced from the joint law
$\pi\colon S\rightarrow X\times Y$.
Different joint laws give different conditional laws, and
properties such as the tower property, Bayesian coherence and Bayesian filters are statements about
how conditional laws of related joint laws fit together, which the traditional notation
$p(y\mid x)$ leaves implicit.  Note that our theory is based on probability measures,
so density functions are not required.

The existence and uniqueness of the conditional law is the parameterized version of the
classical disintegration theorem, see for example \cite[Theorem~1.25]{Kall17}.
The proof including our case is given in \cite{For21}.
We give a self-contained proof in Appendix \ref{sec: appendix A}.

\begin{theorem}\label{thm: CL}
	Consider a joint law $\pi\colon S\rightarrow X\times Y$,
	 where $X$ is countably generated.
	\begin{enumerate}[label=(\roman*)]
		\item If $Y$ is countably generated, and if a conditional law
		 $\pi(Y\mid X)$ of $Y$ given $X$ for $\pi$ exists,
		then such $\pi(Y\mid X)$ is unique $\pi_X(s)$-almost everywhere for any $s\in S$.
		\item In particular, if $Y$ is standard Borel, then
		a conditional law
		 $\pi(Y\mid X)$ of $Y$ given $X$ for $\pi$ exists,
		 and is determined uniquely $\pi_X(s)$-almost everywhere for any $s\in S$.
	\end{enumerate}
\end{theorem}

\begin{proof}
	See Appendix \ref{sec: appendix A}, where the theorem is restated as \cref{thm: main}.
\end{proof}

From now on, we assume that the spaces $X,Y,Z,V,W$ appearing as factors of joint laws are countably generated,
while the parameter space $S$ is an arbitrary measurable space.
This convention is stable under taking finite or countable products and replacing $S$ by $S\times X$.
By \cref{thm: CL}, a conditional law is determined uniquely almost everywhere whenever it exists.
Existence will either be assumed explicitly or guaranteed by a standard Borel hypothesis.
Equalities between conditional laws are understood in the almost everywhere sense, so that
for a stochastic morphism $\kappa\colon S\times X\rightarrow Y$,
the identity
\[
	\kappa=\pi(Y\mid X)
\]
means that $\kappa=\pi(Y\mid X)$ for $\pi_X(s)$-almost every $x\in X$ for any $s\in S$,
and holds for any choice of the conditional laws involved.

\begin{remark}\label{rem: prX}
	For $S=\bsone$, let $\pi$ be the joint law of random variables $\bX$ and $\bY$ on
	a probability space $(\Omega,\sF,\bbP)$ with values in $X$ and $Y$.
	If the conditional law $\pi(Y\mid X)\colon X\rightarrow Y$ exists, then we have
	\[
		\bbP(\bY\in B\mid\bX=x)=\pi(Y\mid X)(x)(B)
	\]
	for $\pi_X$-almost every $x\in X$ for any measurable $B\subset Y$, that is,
	$\pi(Y\mid X)$ is a regular conditional distribution of $\bY$ given $\bX=x$.
	The factor $X$ plays the role of the random variable $\bX$, realized by the
	coordinate projection $\pr_X$, and our notation can be read as in classical probability theory.
\end{remark}

\subsection{Augmentations of Stochastic Morphisms}

In the characterization $\eqref{eq: CL}$ of the conditional law,
we require $(\id_S, \pi_X)$, which holds
information of the input.  We next introduce the augmentation,
which preserves the input of stochastic morphisms.

\begin{definition}
	Given a stochastic morphism
	$
		\kappa\colon S\rightarrow X,
	$
	we define the \emph{augmentation} $\kappa^\wedge$ of $\kappa$ by
	\[
		\kappa^\wedge\coloneqq(\id_S,\kappa)\colon S\rightarrow S\times X,
	\]
	which augments the output by the input.
\end{definition}
Categorically, the augmentation $\kappa^\wedge$ is simply the graph of $\kappa$.
By definition, we have $(\kappa^\wedge)_X=\kappa$, so that $\kappa$ is determined by $\kappa^\wedge$.
We use augmentations systematically in compositions, so that the characterization
\eqref{eq: CL} of the conditional law becomes the strictly commutative diagram
\eqref{eq: CL aug} below.
The augmentation satisfies the following properties.

\begin{lemma}[Properties of Augmentations]\label{lem: augmentations}
	\begin{enumerate}[label=(\alph*)]
		\item\label{aug: P} (Push-forward)
		For stochastic morphisms $\kappa\colon S\rightarrow X$ and $a\colon S\times X\rightarrow Y$, we have
		\[
			(\pr_S,a)\circ\kappa^\wedge=(a\circ\kappa^\wedge)^\wedge\colon S\rightarrow S\times Y.
		\]
	\item\label{aug: M} (Marginals)
	For a stochastic morphism $\kappa\colon S\rightarrow X\times Y$, we have
	\[
		(\kappa^\wedge)_{S\times X}=(\id_S,\kappa)_{S\times X}=(\id_S,\kappa_X)=(\kappa_X)^\wedge.
	\]
	\item\label{aug: I} (Independence)
	For stochastic morphisms $p\colon S\rightarrow X$ and $q\colon S\rightarrow Y$, we have
	\[
		(q\circ\pr_S)^\wedge\circ p^\wedge=(\id_S,p,q)=(p,q)^\wedge,
	\]
	where $q\circ\pr_S\colon S\times X\rightarrow Y$.
	\end{enumerate}
\end{lemma}

\begin{proof}
	Each property is verified by evaluating both sides on measurable rectangles, since a probability measure
	on a product is determined by its values on rectangles.
	For $\ref{aug: P}$, for measurable sets $A\subset S$ and $B\subset Y$, we have
	\begin{align*}
		((\pr_S,a)\circ\kappa^\wedge)(s)(A\times B)
		&=\int_{S\times X}\delta_{s'}(A)\,a(s',x)(B)\,\kappa^\wedge(s)(ds',dx)\\
		&=\int_{X}\delta_{s}(A)\,a(s,x)(B)\,\kappa(s)(dx)\\
		&=\delta_s(A)\,(a\circ\kappa^\wedge)(s)(B)
		=(a\circ\kappa^\wedge)^\wedge(s)(A\times B),
	\end{align*}
	since $\kappa^\wedge(s)(ds',dx)=\delta_s(ds')\kappa(s)(dx)$.
	Similarly for $\ref{aug: M}$ and $\ref{aug: I}$.
\end{proof}

The properties apply to any parameter space, for example to $S\times X$ in place of $S$.
By applying \cref{lem: augmentations}~$\ref{aug: P}$ to \eqref{eq: CL} with $\kappa=\pi_X$ and $a=(\pr_X,\pi(Y\mid X))$, we see that
\begin{align*}
	\pi^\wedge&= ((\pr_X, \pi(Y\mid X))\circ (\id_S,\pi_X))^\wedge
	=((\pr_X, \pi(Y\mid X))\circ (\pi_X)^\wedge)^\wedge\\
	&=(\pr_S,(\pr_X, \pi(Y\mid X)))\circ(\pi_X)^\wedge
	=\pi(Y\mid X)^\wedge\circ(\pi_X)^\wedge\colon S\rightarrow S\times X\times Y.
\end{align*}
Thus noting that a stochastic morphism is determined by its augmentation, \eqref{eq: CL} is equivalent to
\[
	\pi^\wedge=\pi(Y\mid X)^\wedge\circ(\pi_X)^\wedge\colon S\rightarrow S\times X\times Y.
\]
In other words, the diagram
\begin{equation}\label{eq: CL aug}
	\xymatrix{
		S \ar[dd]_{\pi^\wedge} \ar[dr]^{(\pi_X)^\wedge} &  \\
		&  S\times X \ar[dl]^{\pi(Y\mid X)^\wedge}\\
		S\times X\times Y&
	}
\end{equation}
is commutative.
We stress that the diagram commutes strictly as an equality of stochastic morphisms
for any choice of the conditional law, and not merely up to almost-everywhere equivalence.

More generally, let $I$ be a countable index set and let $\{X_i\}_{i\in I}$ be a system of countably generated measurable spaces.
For any $J\subset I$, let
$
	X^J\coloneqq\prod_{i\in J}X_i,
$
which is again countably generated.
Consider a joint law
$
	\pi\colon S\rightarrow X^I.
$
For any $J,K\subset I$ which are disjoint, we let
\[
	\pi(X^J\mid X^K)\colon S\times X^K\rightarrow X^J
\]
be the conditional law of $X^J$ given $X^K$ for the marginal $\pi_{X^K\times X^J}$, whenever it exists.
In particular, $\pi(X^J\mid X^K)$ depends only on the marginal $\pi_{X^K\times X^J}$,
so that the conditional laws of a joint law and of its marginals are \emph{consistent},
and coincide wherever both are defined.
We extend the notation to $J,K\subset I$ which are not necessarily disjoint by copying the shared coordinates, so that
\[
	\pi(X^J\mid X^K)\coloneqq(\pr_{X^{J\cap K}},\,\pi(X^{J\setminus K}\mid X^K))\colon S\times X^K\rightarrow X^{J\cap K}\times X^{J\setminus K}=X^J
\]
up to the order of the factors. 
By definition,
the probability measure $\pi(X^J\mid X^K)(s,x)$ 
for $s\in S$ and $x\in X^K$
is the product of the Dirac measure at the $X^{J\cap K}$-coordinates of $x$
and the measure $\pi(X^{J\setminus K}\mid X^K)(s,x)$.
In particular, for a joint law $\pi\colon S\rightarrow X\times Y$, we have
$\pi(X\times Y\mid X)=(\pr_X,\pi(Y\mid X))\colon S\times X\rightarrow X\times Y$ and $\pi(Y\mid X)^\wedge=(\id_S,\pi(X\times Y\mid X))$,
and the characterization \eqref{eq: CL} of the conditional law reads
\begin{equation}\label{eq: CL copy}
	\pi=\pi(X\times Y\mid X)\circ(\pi_X)^\wedge.
\end{equation}
We also take the convention that a stochastic morphism defined on a subproduct, such as $\kappa\colon S\times W\rightarrow Y$,
is regarded as a stochastic morphism on a larger product, such as $S\times W\times X$, through the projection,
whenever the larger domain is determined by the context, namely by the other side of an equality,
by a composition, or by a product $(p,q)$.
We do not apply this convention inside an augmentation, so that the domain of $\kappa^\wedge$ is always that of $\kappa$.

Given an inclusion of index sets $I\subset I'$
and a system of countably generated measurable spaces $\{X_i\}_{i\in I'}$,
we say that a joint law $\pi'\colon S\rightarrow X^{I'}$ is an \emph{extension}
of $\pi\colon S\rightarrow X^{I}$ if $\pi'_{X^{I}}=\pi$.
Then for $J,K\subset I$, we have
\[
	\pi'(X^J\mid X^K)=\pi(X^J\mid X^K)
\]
in this case.

%
\subsection{Properties of Conditional Laws}
%

In this subsection, we prove fundamental properties of the conditional law.
The conditional law satisfies the following.

\begin{lemma}[Properties of Conditional Laws]\label{lem: properties}
	\begin{enumerate}[label=(\alph*)]
		\item\label{CL: P} (Push-forward)
		Let $\pi\colon S\rightarrow X\times Y$
		be a joint law such that the conditional law  $\pi(Y\mid X)$ exists,
		and consider the stochastic morphism $a\colon Y\rightarrow Z$.
		Then we have
		\[
			a\circ\pi(Y\mid X)=((\id_X\times a)\circ\pi)(Z\mid X).
		\]
		More generally, for $a\colon S\times X\times Y\rightarrow Z$ 
		parameterized by $S\times X$, we have
		\[
			a\circ\pi(Y\mid X)^\wedge=((\pr_X,a)\circ\pi^\wedge)(Z\mid X).
		\]
		In particular, for $a\colon X\times Y\rightarrow Z$ parameterized by $X$, we have
		\[
			a\circ\pi(X\times Y\mid X)=((\pr_X,a)\circ\pi)(Z\mid X).
		\]
		Each case guarantees the existence of the conditional law on the right hand side.
		\item\label{CL: M} (Marginals)
		Let $\pi\colon S\rightarrow X\times Y\times Z$ be a joint law such 
		that $\pi(Y\times Z\mid X)$ exists.
		Then we have
		\[
			\pi(Y\times Z\mid X)_Y=\pi(Y\mid X),
		\]
		where
		$\pi(Y\times Z\mid X)_Y\coloneqq\pr_Y\circ\pi(Y\times Z\mid X)$ is the $Y$-marginal of
		$\pi(Y\times Z\mid X)$.
		\item\label{CL: C} (Chain rule)
		Let $\pi\colon S\rightarrow X\times Y\times Z$ be a joint law such that
		$\pi(Y\mid X)$ and $\pi(Z\mid X\times Y)$ exist.  Then $\pi(Z\mid X)$ exists, and
		\[
			\pi(Z\mid X)=\pi(Z\mid X\times Y)\circ\pi(Y\mid X)^\wedge.
		\]
	\end{enumerate}
\end{lemma}

Classically, the chain rule \ref{CL: C} reads $p(z\mid x)=\int_Y p(z\mid x,y)\,p(y\mid x)\,dy$ for densities.

\begin{proof}
	We first consider the case $a\colon S\times X\times Y\rightarrow Z$.
	 Let $q\coloneqq\pi(Y\mid X)\colon S\times X\rightarrow Y$.
	By \cref{lem: augmentations}~$\ref{aug: P}$ with $S\times X$ as the parameter space, we have
	$(\pr_{S\times X},a)\circ q^\wedge=(a\circ q^\wedge)^\wedge\colon S\times X\rightarrow S\times X\times Z$.  Composing with the projection
	$\pr_{X\times Z}=\pr_X\times\id_Z$ and using \eqref{eq: product}, we obtain
	\[
		(\pr_X,a)\circ q^\wedge=(\pr_X,\,a\circ q^\wedge)\colon S\times X\rightarrow X\times Z.
	\]
	Hence by \eqref{eq: CL aug} for $\pi$, with $q=\pi(Y\mid X)$, we have
	\[
		(\pr_X,\,a\circ q^\wedge)\circ(\pi_X)^\wedge=(\pr_X,a)\circ q^\wedge\circ(\pi_X)^\wedge=(\pr_X,a)\circ\pi^\wedge\colon S\rightarrow X\times Z.
	\]
	Since the $X$-marginal of $(\pr_X,a)\circ\pi^\wedge$ is $\pr_X\circ\pi^\wedge=\pi_X$,
	this is the defining identity \eqref{eq: CL} of the conditional law of $Z$ given $X$
	for $(\pr_X,a)\circ\pi^\wedge$, which proves the assertion for $a\colon S\times X\times Y\rightarrow Z$.
	For $a\colon X\times Y\rightarrow Z$, apply this to $a\circ\pr_{X\times Y}$, noting that
	$(\pr_X,a\circ\pr_{X\times Y})\circ\pi^\wedge=(\pr_X,a)\circ\pi\colon S\rightarrow X\times Z$ and
	$\pr_{X\times Y}\circ q^\wedge=(\pr_X,q)\colon S\times X\rightarrow X\times Y$ by \eqref{eq: product}.
	The case $a\colon Y\rightarrow Z$ follows by replacing $a$ with $a\circ\pr_Y$.

	We next prove $\ref{CL: M}$.  Let $\kappa=\pi(Y\times Z\mid X)$.
	By \eqref{eq: CL aug} for $\pi$ and by \cref{lem: augmentations}~$\ref{aug: M}$ applied twice, we have
	\[
		(\pi_{X\times Y})^\wedge=(\pi^\wedge)_{S\times X\times Y}
		=\pr_{S\times X\times Y}\circ\kappa^\wedge\circ(\pi_X)^\wedge
		=(\kappa^\wedge)_{S\times X\times Y}\circ(\pi_X)^\wedge
		=(\kappa_Y)^\wedge\circ(\pi_X)^\wedge,
	\]
	which is the commutative diagram \eqref{eq: CL aug} of the conditional law of $Y$ given $X$ for $\pi_{X\times Y}$.
	The assertion follows from the uniqueness of the conditional law.

	For $\ref{CL: C}$, the defining identity \eqref{eq: CL} of $\pi(Z\mid X\times Y)$ reads
	$\pi=(\pr_{X\times Y},\pi(Z\mid X\times Y))\circ(\pi_{X\times Y})^\wedge$.
	Composing with $\pr_{X\times Z}$ and using \eqref{eq: product}, we obtain
	\[
		\pi_{X\times Z}=(\pr_X,\pi(Z\mid X\times Y))\circ(\pi_{X\times Y})^\wedge\colon S\rightarrow X\times Z,
	\]
	which is the push-forward of $\pi_{X\times Y}$ by $a=\pi(Z\mid X\times Y)\colon S\times X\times Y\rightarrow Z$.
	Hence $\ref{CL: P}$ applied to $\pi_{X\times Y}$ gives
	$\pi(Z\mid X)=\pi(Z\mid X\times Y)\circ\pi(Y\mid X)^\wedge$.
\end{proof}

\begin{theorem}[The Tower Property]\label{thm: tower}
	Consider a joint law $\pi\colon S\rightarrow X\times Y\times Z$.
	If $\pi(Y\times Z\mid X)\colon S\times X\rightarrow Y\times Z$ is the
	conditional law of $Y\times Z$ given $X$ for $\pi$, and
	\[
		\pi(Y\times Z\mid X)(Z\mid Y)\colon S\times X\times Y\rightarrow Z
	\]
	is the conditional law of $Z$
	given $Y$ for $\pi(Y\times Z\mid X)$, with parameter space $S\times X$,
	then the conditional law $\pi(Z\mid X\times Y)$ of $Z$ given $X\times Y$ for $\pi$ exists, and
	we have
	\[
		\pi(Y\times Z\mid X)(Z\mid Y)=\pi(Z\mid X\times Y)
	\]
	for $\pi_{X\times Y}(s)$-almost every $(x,y)\in X\times Y$ for any $s\in S$.
\end{theorem}

\begin{proof}
	Let $\kappa\coloneqq\pi(Y\times Z\mid X)\colon S\times X\rightarrow Y\times Z$, so that the left hand side is
	$\kappa(Z\mid Y)$.
	We have a commutative diagram
	\[
		\xymatrix@C=20pt@R=20pt{
			  S\ar[ddd]_{\pi^\wedge} \ar[dr]_{(\pi_X)^\wedge} \ar@/^6ex/[ddrr]^{(\pi_{X\times Y})^\wedge} &  &\\
			  &    S\times X\ar[ddl]_{\kappa^\wedge}  \ar[dr]_{(\kappa_Y)^\wedge}&   \\
			 &  & S\times X \times Y\ar[lld]^{\kappa(Z\mid Y)^\wedge}\\
			  S\times X\times Y \times Z. &  &
		}
	\]
	Here, the commutativity of the bottom left and right triangles is the commutative diagram
	\eqref{eq: CL aug} of the conditional laws $\kappa$ of $\pi$ and $\kappa(Z\mid Y)$ of $\kappa$,
	and that of the upper triangle is the commutative diagram \eqref{eq: CL aug} of the conditional law
	$\pi(Y\mid X)=\kappa_Y$ for $\pi_{X\times Y}$, by \cref{lem: properties}~$\ref{CL: M}$.
	Hence
	\[
		\pi^\wedge=\kappa^\wedge\circ(\pi_X)^\wedge=\kappa(Z\mid Y)^\wedge\circ(\kappa_Y)^\wedge\circ(\pi_X)^\wedge
		=\kappa(Z\mid Y)^\wedge\circ(\pi_{X\times Y})^\wedge,
	\]
	which is the commutative diagram \eqref{eq: CL aug} of the conditional law of $Z$ given $X\times Y$ for $\pi$.
	Hence $\kappa(Z\mid Y)$ is a conditional law of $Z$ given $X\times Y$ for $\pi$, which proves the existence,
	and the assertion follows from the uniqueness of the conditional law.
\end{proof}

We next formulate conditional independence in terms of conditional laws.

\begin{definition}\label{def: CI}
	Let $\pi\colon S\rightarrow W\times X\times Y$ be a joint law.
	We say that $X$ and $Y$ are \emph{conditionally independent given $W$} for $\pi$, and write
	\[
		X\perp_\pi Y\mid W,
	\]
	if the conditional laws $\pi(X\mid W)$ and $\pi(Y\mid W)$ exist and their product
	$(\pi(X\mid W),\pi(Y\mid W))\colon S\times W\rightarrow X\times Y$
	gives a conditional law of $X\times Y$ given $W$ for $\pi$, that is,
	\[
		\pi(X\times Y\mid W)=(\pi(X\mid W),\pi(Y\mid W)).
	\]
\end{definition}

Classically, this reads $p(x,y\mid w)=p(x\mid w)\,p(y\mid w)$ for densities.
The condition does not depend on the choice of $\pi(X\mid W)$ and $\pi(Y\mid W)$,
since the right hand side of \eqref{eq: CL} does not change if the conditional law is modified on $\pi_W(s)$-null sets.
If $W=\bsone$, the condition reads $\pi=(\pi_X,\pi_Y)$, which is the independence of $X$ and $Y$ for $\pi$.
The following lemma expresses conditional independence as the irrelevance of $X$ for the conditional law of $Y$.

\begin{lemma}[Conditional Independence]\label{lem: CI}
	Let $\pi\colon S\rightarrow W\times X\times Y$ be a joint law.
	Then the following are equivalent.
	\begin{enumerate}[label=(\roman*)]
		\item\label{CI: product} $X\perp_\pi Y\mid W$.
		\item\label{CI: drop} The conditional laws $\pi(X\mid W)$, $\pi(Y\mid W)$ and $\pi(Y\mid W\times X)$ exist, and
		\[
			\pi(Y\mid W\times X)=\pi(Y\mid W).
		\]
	\end{enumerate}
	By symmetry, these are also equivalent to the existence of $\pi(X\mid W)$, $\pi(Y\mid W)$ and $\pi(X\mid W\times Y)$
	with $\pi(X\mid W\times Y)=\pi(X\mid W)$.
\end{lemma}

\begin{proof}
	In both cases $\pi(X\mid W)$ and $\pi(Y\mid W)$ exist.
	By \cref{lem: augmentations}~$\ref{aug: I}$ with $S\times W$ as the parameter space,
	applied to $p=\pi(X\mid W)$ and $q=\pi(Y\mid W)$, we have
	\[
		(\pi(X\mid W),\pi(Y\mid W))^\wedge=(\pi(Y\mid W)\circ\pr_{S\times W})^\wedge\circ\pi(X\mid W)^\wedge.
	\]
	Here $\pi(Y\mid W)\circ\pr_{S\times W}\colon S\times W\times X\rightarrow Y$.
	Composing with $(\pi_W)^\wedge$ and using the commutative diagram \eqref{eq: CL aug}
	 of $\pi(X\mid W)$
	giving
	$(\pi_{W\times X})^\wedge=\pi(X\mid W)^\wedge\circ(\pi_W)^\wedge$, we obtain
	\begin{equation}\label{eq: CI}
		(\pi(X\mid W),\pi(Y\mid W))^\wedge\circ(\pi_W)^\wedge=(\pi(Y\mid W)\circ\pr_{S\times W})^\wedge\circ(\pi_{W\times X})^\wedge.
	\end{equation}
	By \eqref{eq: CL aug} and the uniqueness of the conditional law,
	\ref{CI: product} holds if and only if the left hand side of \eqref{eq: CI} is $\pi^\wedge$,
	and \ref{CI: drop} holds if and only if the right hand side of \eqref{eq: CI} is $\pi^\wedge$.
	This proves the equivalence.
\end{proof}

Conditional independence arises when a variable depends solely on the conditioning variable.

\begin{corollary}\label{cor: CI}
	Let $\pi\colon S\rightarrow W\times X$ be a joint law such that $\pi(X\mid W)$ exists,
	let $b\colon W\rightarrow Y$ be a stochastic morphism, and let
	\[
		\pi'\coloneqq ((\id_W,b)\times \id_X)\circ\pi
	\]
	be the extension of $\pi$ in which $Y$ is generated from $W$ by $b$.
	Then $X\perp_{\pi'} Y\mid W$, and hence
	\[
		\pi'(X\mid W\times Y)=\pi(X\mid W).
	\]
\end{corollary}

\begin{proof}
	Since $(\id_W,b)\times\id_X=(\pr_W,\,b\times\id_X)$, \cref{lem: properties}~$\ref{CL: P}$ applied to $\pi$,
	with $W$ and $X$ in place of $X$ and $Y$, $Z\coloneqq Y\times X$ and
	$a\coloneqq b\times\id_X\colon W\times X\rightarrow Y\times X$, together with \eqref{eq: product}, gives
	\[
		\pi'(Y\times X\mid W)=(b\times\id_X)\circ(\pr_W,\pi(X\mid W))=(b,\,\pi(X\mid W))\colon S\times W\rightarrow Y\times X.
	\]
	By \cref{lem: properties}~$\ref{CL: M}$, we have $\pi'(Y\mid W)=b$ and $\pi'(X\mid W)=\pi(X\mid W)$,
	so that $\pi'(X\times Y\mid W)=(\pi'(X\mid W),\pi'(Y\mid W))$, which is $X\perp_{\pi'}Y\mid W$.
	The last assertion follows from \cref{lem: CI}.
\end{proof}

%
%
%
\section{The Bayesian inverse}\label{sec: BI}
%
%
%

In this section, we first define the Bayesian inverse, and then prove Bayesian coherence,
which states that updating observations successively is the same as updating them at once.

%
\subsection{Bayesian Inverse}
%

We define the Bayesian inverse as follows.

\begin{definition}\label{def: BI}
	Consider stochastic morphisms $\mu\colon S\rightarrow X$ and $h\colon X\rightarrow Y$.
	The \emph{Bayesian inverse} of $h$ with respect to prior $\mu$ is a stochastic morphism
	\[
		h^\mu\colon S\times Y\rightarrow X
	\]
	characterized by the equality
	\begin{equation}\label{eq: BI}
			(h,\id_X)\circ\mu=(\pr_Y,h^\mu)\circ (\id_S, h\circ \mu)
	\end{equation}
	as morphisms from $S$ to $Y\times X$.
\end{definition}

This is essentially the Bayesian inverse defined by Cho and Jacobs in \cite[\S3]{CJ19},
generalized to allow for parameterized priors.  
The Bayesian inverse can be constructed as a conditional law as follows.

\begin{proposition}\label{prop: construct}
	Suppose $X$ is standard Borel.
	Consider stochastic morphisms $\mu\colon S\rightarrow X$ and $h\colon X\rightarrow Y$.
	Then for a joint law $\pi$ parameterized by $S$ such that $\pi_{Y\times X}=(h,\id_X)\circ\mu$,
	we have
	\[
		h^\mu=\pi(X\mid Y)
	\]
	for $\pi_Y(s)$-almost every $y\in Y$ for any $s\in S$.
\end{proposition}

\begin{proof}
	Since $\pi(X\mid Y)$ is the conditional law of $X$ given $Y$ for the marginal $\pi_{Y\times X}$,
	we may replace $\pi$ by $\pi_{Y\times X}=(h,\id_X)\circ\mu\colon S\rightarrow Y\times X$.
	By our assumption, there exists
	a conditional law $\pi(X\mid Y)\colon S\times Y\rightarrow X$.
	Then \eqref{eq: BI} is simply the characterization
	\[
		\pi=(\pr_Y, \pi(X\mid Y))\circ (\id_S,\pi_Y)
	\]
	of \eqref{eq: CL} for $\pi(X\mid Y)$ noting that $\pi_Y=h\circ\mu$.
	Our assertion follows from the uniqueness of the conditional law.
\end{proof}

If $S=\bsone$, then $\mu$ simply gives a probability measure on $X$, and $h^\mu\colon Y\to X$ is the Bayesian inverse of Cho--Jacobs \cite{CJ19} for $\mu$ along $h$.
The Bayesian inverse in the measure-theoretic setting is also treated in \cite{CS14},
and its relation to disintegrations in general Markov categories in \cite{P20}.
In \cite[Remark 13.10]{F20},
the notation $h^\dagger$ is used for the Bayesian inverse of $h$ due to its interpretation
as a symmetric monoidal dagger functor.  We use $h^\mu$
to identify the prior, which is important in a recursive setting (see \cref{thm: BC}).

%
\subsection{Bayesian Coherence}
%

In this subsection, in order to illustrate the workings of the Bayesian inverse,
we prove Bayesian coherence, which asserts that the posterior given two observations
can be calculated iteratively, the posterior given the first observation serving as the prior for the second.
Let $\mu\colon S\rightarrow W$ be a prior, and consider observations $g\colon W\rightarrow X$ and $h\colon W\rightarrow Y$,
where $W$ is standard Borel.
\[
	\xymatrix@C=50pt@R=10pt{
			&   X  \\
			W  \ar[ur]^g \ar[dr]_h&\\
			&  Y
	}
\]
We define the \emph{batch posterior}
\[
	(g,h)^\mu\colon S\times X\times Y\rightarrow W
\]
to be the Bayesian inverse of $(g,h)\colon W\rightarrow X\times Y$ with respect
to the prior $\mu\colon S\rightarrow W$.
Bayesian coherence asserts that the batch posterior can be calculated iteratively.
Compositionality of Bayesian inversion is also studied in the language of lenses \cite{S20,BHS23},
where the Bayesian inverse is indexed by the prior.
Let
\[
	g^\mu\colon S\times X\rightarrow W
\]
be the Bayesian inverse of $g\colon W\rightarrow X$ with respect to the prior $\mu\colon S\rightarrow W$,
and
\[
	h^{g^\mu}\colon S\times X\times Y\rightarrow W
\]
be the Bayesian inverse of $h$ with respect to the prior $g^\mu\colon S\times X\rightarrow W$,
with parameter space $S\times X$.

\begin{theorem}[Bayesian Coherence]\label{thm: BC}
	Let the notations be as above.  Then we have
	\[
		(g,h)^\mu=h^{g^\mu}
	\]
	for $((g,h)\circ\mu)(s)$-almost every $(x,y)\in X\times Y$ for any $s\in S$.
\end{theorem}

The theorem asserts that the posterior given both observations is obtained
as the Bayesian inverse of $h$, by taking the prior to be the posterior
$g^\mu$ for the first observation.
The parameterized prior in \cref{def: BI} allows this iteration,
taking the previous posterior as the prior.

\begin{proof}
	Let  $\pi\coloneqq (g,\id_W)\circ\mu$ and consider the extension
	\[
		\pi'=(\id_X\times (h,\id_W))\circ\pi \colon S\rightarrow X\times Y\times W
	\]
	which satisfies $\pi'=((g,h),\id_W)\circ\mu$.
	By \cref{prop: construct}, we have
	$(g,h)^\mu=\pi'(W\mid X\times Y)$ and $g^\mu=\pi(W\mid X)$.
	By \cref{lem: properties}~$\ref{CL: P}$ for $a=(h,\id_W)$, we have
	\[
		\pi'(Y\times W\mid X)=(h,\id_W)\circ g^\mu.
	\]
	By \cref{prop: construct} with $S\times X$ as the parameter space and $g^\mu$ as the prior,
	$h^{g^\mu}$ is the conditional law of $W$ given $Y$ for $(h,\id_W)\circ g^\mu=\pi'(Y\times W\mid X)$, that is,
	\[
		h^{g^\mu}=\pi'(Y\times W\mid X)(W\mid Y).
	\]
	The tower property of \cref{thm: tower} gives
	\[
		\pi'(Y\times W\mid X)(W\mid Y)=\pi'(W\mid X\times Y)=(g,h)^\mu,
	\]
	which proves $h^{g^\mu}=(g,h)^\mu$ as desired.
\end{proof}

\begin{remark}\label{rem: BC}
	The proof uses the joint law $\pi'\colon S\rightarrow X\times Y\times W$ only through the properties
	$\pi'(Y\mid W)=h$ and $X\perp_{\pi'}Y\mid W$,
	which hold by \cref{lem: properties}~$\ref{CL: M}$ and \cref{cor: CI}, since $X$ and $Y$ are generated from $W$ by $g$ and $h$.
	Indeed, for any joint law $\pi\colon S\rightarrow W\times X\times Y$ with $W$ standard Borel
	such that $\pi(Y\mid W)=h$ and $X\perp_\pi Y\mid W$, we have
	\[
		\pi(W\mid X\times Y)=h^{\pi(W\mid X)},
	\]
	the Bayesian inverse of $h$ with respect to the prior $\pi(W\mid X)$ with parameter space $S\times X$.
	By \cref{lem: CI}, we have 
	$\pi(Y\mid W\times X)=h$, so that the commutative diagrams \eqref{eq: CL aug}
	of $\pi(Y\mid W\times X)$ and $\pi(W\mid X)$ together with \cref{lem: augmentations}~$\ref{aug: P}$ give
	$\pi(W\times Y\mid X)=(\id_W,h)\circ\pi(W\mid X)$, and the argument above applies.
	This is Bayesian coherence for a first observation $X$ which may depend on the parameter.
\end{remark}

We extend \cref{thm: BC} to a sequence of observations.
Let $X$ be a standard Borel space, and we fix a prior $\mu\colon S\rightarrow X$.
Let $\{Y_i\}_{i=1}^n$ be a family of measurable
spaces, with stochastic morphisms $h_i\colon X\rightarrow Y_i$ for $i=1,\ldots,n$.
Let
\[
	(h_1,\ldots,h_k)\colon X\rightarrow Y^{[k]}\coloneqq\prod_{i=1}^k Y_i,\qquad 1\leq k\leq n.
\]
We define the batch posterior to be the Bayesian inverse
\[
	(h_1,\ldots,h_k)^\mu\colon S\times Y^{[k]}\rightarrow X
\]
of $(h_1,\ldots,h_k)$ with respect to the prior $\mu$,
which by \cref{prop: construct} is the conditional law $\pi_k(X\mid Y^{[k]})$ of $X$ given $Y^{[k]}$ for
$\pi_k\coloneqq((h_1,\ldots,h_k),\id_X)\circ\mu\colon S\rightarrow Y^{[k]}\times X$.

\begin{corollary}[Bayesian Coherence]\label{cor: BC}
	Let the notations be as above, and let $1\leq k<n$.
	Then the Bayesian inverse
	\[
		h_{k+1}^{(h_1,\ldots,h_k)^\mu}\colon S\times Y^{[k+1]}\rightarrow X
	\]
	of $h_{k+1}$ with respect to the prior $(h_1,\ldots,h_k)^\mu$, with parameter space $S\times Y^{[k]}$,
	is the conditional law $\pi_{k+1}(X\mid Y^{[k+1]})$ of $X$ given $Y^{[k+1]}$ for $\pi_{k+1}$.  In particular, we have
	\[
		(h_1,\ldots,h_{k+1})^\mu=h_{k+1}^{(h_1,\ldots,h_k)^\mu}
	\]
	for $((h_1,\ldots,h_{k+1})\circ\mu)(s)$-almost every $(y_1,\ldots,y_{k+1})\in Y^{[k+1]}$ for any $s\in S$.
\end{corollary}

\begin{proof}
	We apply \cref{thm: BC} with $X$ in place of $W$, with $Y^{[k]}$ and $Y_{k+1}$ in place of $X$ and $Y$,
	and with $g=(h_1,\ldots,h_k)$ and $h=h_{k+1}$.
	Since $Y^{[k+1]}=Y^{[k]}\times Y_{k+1}$ and $(g,h)=(h_1,\ldots,h_{k+1})$,
	\cref{thm: BC} gives
	$(h_1,\ldots,h_{k+1})^\mu=h_{k+1}^{(h_1,\ldots,h_k)^\mu}$,
	and $(h_1,\ldots,h_{k+1})^\mu=\pi_{k+1}(X\mid Y^{[k+1]})$ by \cref{prop: construct}.
\end{proof}

\Cref{cor: BC} is the Bayesian filter for a hidden state which does not move.
The posterior after $k+1$ observations is a single Bayesian inverse,
with the past observations entering only through the prior $(h_1,\ldots,h_k)^\mu$ and the parameter space $S\times Y^{[k]}$.
In \S\ref{sec: filter}, the hidden state evolves between the observations, and the prior is replaced by the predicted prior.

%
%
%
\section{Filtering, Prediction and Smoothing}\label{sec: filter}
%
%
%

In this section, we derive the Bayesian filter, prediction and smoothing for a hidden Markov model with feedback,
whose transitions and observations may depend on the past observations.

\subsection{Hidden Markov Models with Feedback}

In this subsection, we define the hidden Markov model with feedback and its joint law,
and extend the Bayesian inverse to priors parameterized by the past observations.
For measurable spaces $\{Y_i\}_{i\geq1}$ and integers $0\leq k\leq j$, we let
\[
	Y^{(k,j]}\coloneqq\prod_{i=k+1}^jY_i,\qquad Y^{[j]}\coloneqq Y^{(0,j]},
\]
with the convention that $Y^{(k,k]}\coloneqq\bsone$, and similarly for $X^{[j]}$.

\begin{definition}[Hidden Markov Model with Feedback]\label{def: HMP}
	The model consists of a system of standard Borel spaces $\{X_i\}_{i=1}^n$
	with a prior $\mu\colon S\rightarrow X_1$,
	and a family of measurable spaces $\{Y_i\}_{i=1}^n$,
	together with
	stochastic morphisms
	\[
		a_i\colon Y^{[i]}\times X_i\rightarrow X_{i+1},\qquad 1\leq i < n
	\]
	governing the transitions of the hidden state, and
	stochastic morphisms
	\[
		h_i\colon Y^{[i-1]}\times X_i\rightarrow Y_i,\qquad 1\leq i\leq n
	\]
	expressing the observations.
\end{definition}

Our model is defined so that the transitions $a_i$ and observations $h_i$ may depend on past observations,
which we refer to as \emph{feedback}.
Note that $a_i$ may depend on the current observation $y_i$, for instance a control input chosen after observing $y_i$.
The case without feedback, where $a_i$ and $h_i$ factor through $X_i$,
is the usual hidden Markov model expressed by the diagram
\[
	\xymatrix{
			S\ar[r]^\mu &  X_1\ar[d]_{h_1}\ar[r]^{a_1} & X_2 \ar[d]_{h_2} \ar[r]^{a_2}
			& X_3 \ar[d]_{h_3} \ar[r]^{a_3}
			&  \cdots \ar[r]^{a_{n-1}}&  X_n \ar[d]_{h_n} \\
			 &  Y_1 & Y_2 & Y_3 & \cdots &  Y_n.
	}
\]
See \cite{F25} for a categorical framework of the theory in the case without feedback.

For $1\leq k<n$, let
\[
	b_k\coloneqq(h_{k+1},\pr_{X_{k+1}})\circ(\pr_{Y^{[k]}},a_k)\colon Y^{[k]}\times X_k\rightarrow Y_{k+1}\times X_{k+1},
\]
which expresses the transition of the hidden state followed by the observation.
The joint law of the model is the stochastic morphism
\[
	\bspi\colon S\rightarrow Y^{[n]}\times X^{[n]}
\]
defined as $\bspi\coloneqq\bspi_n$, where $\bspi_k\colon S\rightarrow Y^{[k]}\times X^{[k]}$ is defined
iteratively by $\bspi_1\coloneqq (h_1,\id_{X_1})\circ\mu$ and
\[
	\bspi_{k+1}\coloneqq(\pr_{Y^{[k]}},\,b_k,\,\pr_{X^{[k]}})\circ \bspi_{k}
\]
for $1\leq k<n$, where we view $X^{[k+1]}=X_{k+1}\times X^{[k]}$.
All the conditional laws of the model are those of $\bspi$,
with coordinates $Y_1,\ldots,Y_n,X_1,\ldots,X_n$.
Since $\bspi$ is an extension of each of its marginals, they may be computed from any marginal of $\bspi$
containing the relevant coordinates.

The joint law is characterized by its conditional laws as follows.
Here $a_i$ and $h_i$ are regarded as stochastic morphisms on $S\times Y^{[i]}\times X^{[i]}$
and on $S\times Y^{[i-1]}\times X^{[i]}$ through the projections.

\begin{proposition}[Characterization of the Joint Law]\label{prop: characterization}
	The joint law $\bspi$ is the unique joint law $\pi\colon S\rightarrow Y^{[n]}\times X^{[n]}$ such that
	$\pi_{X_1}=\mu$ and
	\begin{align}
		\pi(X_{i+1}\mid Y^{[i]}\times X^{[i]})&=a_i,\quad  1\leq i< n, &
		\pi(Y_i\mid Y^{[i-1]}\times X^{[i]})&=h_i,\quad 1\leq i\leq n,\label{eq: Markov h}
	\end{align}
	in the almost everywhere sense of \S\ref{sec: CL}.
\end{proposition}

\begin{proof}
	We first note that for a prior $\nu\colon S\rightarrow X$ and a stochastic morphism $\kappa\colon X\rightarrow Y$,
	the joint law $\varpi\coloneqq(\id_X,\kappa)\circ\nu\colon S\rightarrow X\times Y$ satisfies $\varpi_X=\nu$ and
	$\varpi=(\pr_X,\kappa)\circ(\id_S,\nu)$ by \eqref{eq: product},
	which is the characterization \eqref{eq: CL} of $\kappa$, regarded as $S\times X\rightarrow Y$,
	as a conditional law of $Y$ given $X$ for $\varpi$.
	Conversely, a joint law $\varpi'\colon S\rightarrow X\times Y$ with $\varpi'_X=\nu$ and $\varpi'(Y\mid X)=\kappa$
	satisfies $\varpi'=(\id_X,\kappa)\circ\nu=\varpi$ by \eqref{eq: CL}, since the right hand side of \eqref{eq: CL}
	does not change if the conditional law is modified on $\varpi'_X(s)$-null sets.

	Let $\varpi_i\coloneqq\bspi_{Y^{[i-1]}\times X^{[i]}}$ for $1\leq i\leq n$, so that $\varpi_1=\mu$.
	By the definition of $b_{i-1}$ and \eqref{eq: product}, the construction of $\bspi_i$ from $\bspi_{i-1}$ factors as
	\[
		\varpi_i=(\id_{Y^{[i-1]}\times X^{[i-1]}},\,a_{i-1})\circ\bspi_{i-1},\qquad
		\bspi_i=(\id_{Y^{[i-1]}\times X^{[i]}},\,h_i)\circ\varpi_i,
	\]
	up to the order of the factors, for $1<i\leq n$ and for $1\leq i\leq n$ respectively.
	Since the conditional laws of $\bspi$ are those of its marginals, the first observation shows that
	$\bspi$ satisfies \eqref{eq: Markov h}, and $\bspi_{X_1}=\mu$ by construction.

	Conversely, let $\pi$ be a joint law such that $\pi_{X_1}=\mu$ satisfying \eqref{eq: Markov h}.
	We prove $\pi_{Y^{[i]}\times X^{[i]}}=\bspi_i$ by induction on $i$, which for $i=n$ gives $\pi=\bspi$.
	Suppose $\pi_{Y^{[i-1]}\times X^{[i-1]}}=\bspi_{i-1}$ if $i>1$.
	Then the second observation applied to \eqref{eq: Markov h} for $a$ and $i-1$ gives $\pi_{Y^{[i-1]}\times X^{[i]}}=\varpi_i$,
	where for $i=1$ this is the assumption $\pi_{X_1}=\mu=\varpi_1$,
	and then applied to \eqref{eq: Markov h} for $h$ gives $\pi_{Y^{[i]}\times X^{[i]}}=\bspi_i$.
\end{proof}

\begin{remark}
	\Cref{prop: characterization} characterizes the model without reference to the construction of 
	$\bspi$.
	Given the whole past, the next state and the next observation are generated by $a_i$ and $h_i$, so that
	the conditional law of $X_{i+1}$ given $Y^{[i]}\times X^{[i]}$ depends on the past 
	only through $Y^{[i]}\times X_i$,
	and that of $Y_i$ given $Y^{[i-1]}\times X^{[i]}$ only through $Y^{[i-1]}\times X_i$.
	This is the Markov property of the model with feedback.
	In the case without feedback, \eqref{eq: Markov h} says that these conditional laws
	depend on the past only through $X_i$, which are the conditional independences characterizing
	hidden Markov models in \cite[Proposition~3.7]{F25}.
	With feedback, these conditional independences
	given $X_i$ alone fail, while \eqref{eq: Markov h} continues to hold.
\end{remark}

We consider the problem of computing the law of the hidden state $X_j$ given observations $Y^{[k]}$
for $\bspi$.
The \emph{filter} is the problem of computing the \emph{posterior}
$\bspi(X_k\mid Y^{[k]})\colon S\times Y^{[k]}\rightarrow X_k$ of the hidden state given the
observations iteratively in $k$. More generally, $\bspi(X_j\mid Y^{[k]})$ for $j>k$ is the \emph{prediction} and
for $j<k$ the \emph{smoothing} of the hidden state.

We write $\pi_k\coloneqq\bspi_{Y^{[k]}\times X_k}$ for the marginal, so that $\pi_k(X_k\mid Y^{[k]})=\bspi(X_k\mid Y^{[k]})$.
By construction of $\bspi$ and \eqref{eq: product}, we have $\pi_1=(h_1,\id_{X_1})\circ\mu$,
\begin{equation}\label{eq: recursion}
	\pi_{k+1}=(\pr_{Y^{[k]}},b_k)\circ \pi_{k},\qquad 1\leq k<n.
\end{equation}

Since the transitions and observations may depend on the past observations,
we need to generalize the Bayesian inverse defined in \cref{def: BI} to stochastic morphisms parameterized by a measurable space.

\begin{definition}\label{def: BI V}
	Let $V$ be a measurable space, and
	consider stochastic morphisms $\mu\colon S\times V\rightarrow X$ and $h\colon V\times X\rightarrow Y$.
	The Bayesian inverse of $h$ with respect to prior $\mu$, with parameter space $S\times V$,
	is a stochastic morphism
	\[
		h^\mu\colon S\times V\times Y\rightarrow X
	\]
	characterized by the equality
	\begin{equation}\label{eq: BI V}
			(h,\pr_X)\circ\mu^\wedge=(\pr_Y,h^\mu)\circ (h\circ\mu^\wedge)^\wedge
	\end{equation}
	as morphisms from $S\times V$ to $Y\times X$.
\end{definition}

If $V=\bsone$, then $\mu^\wedge=(\id_S,\mu)$ and $h$ is regarded as a stochastic morphism on $S\times X$ through $\pr_X$,
so that by \eqref{eq: product} this is \cref{def: BI}.
The parameterized Bayesian inverse can be constructed as a conditional law, as follows.

\begin{proposition}\label{prop: construct V}
	Suppose $X$ is standard Borel.
	Let $V$ be a measurable space, and
	consider stochastic morphisms $\mu\colon S\times V\rightarrow X$ and $h\colon V\times X\rightarrow Y$.
	Then for a joint law $\pi$ parameterized by 
	$S\times V$ such that $\pi_{Y\times X}=(h,\pr_X)\circ\mu^\wedge$,
	we have
	\[
		h^\mu=\pi(X\mid Y)
	\]
	for $\pi_Y(s,v)$-almost every $y\in Y$ for any $(s,v)\in S\times V$.
\end{proposition}

\begin{proof}
	Since $\pi(X\mid Y)$ is the conditional law of $X$ given $Y$ for the marginal $\pi_{Y\times X}$,
	we may replace $\pi$ by $\pi_{Y\times X}=(h,\pr_X)\circ\mu^\wedge\colon S\times V\rightarrow Y\times X$.
	By our assumption, there exists
	a conditional law $\pi(X\mid Y)\colon S\times V\times Y\rightarrow X$, with parameter space $S\times V$.
	Then \eqref{eq: BI V} is simply the characterization
	\[
		\pi=(\pr_Y, \pi(X\mid Y))\circ (\id_{S\times V},\pi_Y)
	\]
	of \eqref{eq: CL} for $\pi(X\mid Y)$ noting that $\pi_Y=h\circ\mu^\wedge$.
	Our assertion follows from the uniqueness of the conditional law.
\end{proof}

%
\subsection{Bayesian Filtering}
%

In this subsection,
we calculate the conditional law of the filter.
For $1\leq k<n$, we let
\[
	\mu_{k+1}\coloneqq a_k\circ\bspi(Y^{[k]}\times X_k\mid Y^{[k]})\colon S\times Y^{[k]}\rightarrow X_{k+1}
\]
be the \emph{predicted prior}, obtained from
the posterior at time $k$ by applying the transition $a_k$.
Here $\bspi(Y^{[k]}\times X_k\mid Y^{[k]})=(\pr_{Y^{[k]}},\bspi(X_k\mid Y^{[k]}))$ carries the past observations to $a_k$.

\begin{theorem}[Bayesian Filter]\label{thm: BF}
	Let the notations be as above.  Then we have
	\[
		\bspi(X_{k+1}\mid Y^{[k+1]})= h_{k+1}^{\mu_{k+1}},\qquad 1\leq k<n,
	\]
	for $\bspi_{Y^{[k+1]}}(s)$-almost every $(y_1,\ldots,y_{k+1})\in Y^{[k+1]}$
	for any $s\in S$,
	where $h_{k+1}^{\mu_{k+1}}$ is the Bayesian inverse of $h_{k+1}$
	for the prior $\mu_{k+1}$, with $V=Y^{[k]}$ in \cref{def: BI V}.
\end{theorem}

\begin{proof}
	Let $1\leq k<n$, and let
	$\kappa\coloneqq b_k\circ\bspi(Y^{[k]}\times X_k\mid Y^{[k]})
	\colon S\times Y^{[k]}\rightarrow Y_{k+1}\times X_{k+1}$.
	By \eqref{eq: recursion} and \cref{lem: properties}~$\ref{CL: P}$ applied to $\pi_k$ and $a=b_k$, we have
	\begin{equation}\label{eq: BF1}
		\kappa=\pi_{k+1}(Y_{k+1}\times X_{k+1}\mid Y^{[k]})=\bspi(Y_{k+1}\times X_{k+1}\mid Y^{[k]}).
	\end{equation}
	Here, the second equality holds since $\bspi$ is an extension of $\pi_{k+1}$.
	On the other hand, by definition of $b_k$ and $\mu_{k+1}$, we have
	$
		\kappa=(h_{k+1},\pr_{X_{k+1}})\circ\mu_{k+1}^\wedge.
	$
	Hence by \cref{prop: construct V} with $V=Y^{[k]}$, we see that
	\begin{equation}\label{eq: BF2}
		h_{k+1}^{\mu_{k+1}}=\kappa(X_{k+1}\mid Y_{k+1}).
	\end{equation}
	Since $\kappa=\bspi(Y_{k+1}\times X_{k+1}\mid Y^{[k]})$ by \eqref{eq: BF1},
	the tower property \cref{thm: tower} and \eqref{eq: BF2} give
	\[
		\bspi(X_{k+1}\mid Y^{[k]}\times Y_{k+1})
		=\kappa(X_{k+1}\mid Y_{k+1})=h_{k+1}^{\mu_{k+1}},
	\]
	which is our assertion since $Y^{[k]}\times Y_{k+1}=Y^{[k+1]}$.
\end{proof}

In the special case that $X_k=X$ and $a_k=\id_X$ for all $k$, and the observations $h_k\colon X\rightarrow Y_k$ do not depend on the past observations,
then $\mu_{k+1}$ is the posterior at time $k$, and \cref{thm: BF}  is exactly \cref{cor: BC}.

%
\subsection{Prediction and Smoothing}
%

Since the transitions may depend on the observations,
the prediction of the hidden state $X_j$ at a future time $j>k$ requires the prediction of the
future observations $Y^{(k,j]}$ as well, which is given by the conditional law
$\bspi(Y^{(k,j]}\times X_j\mid Y^{[k]})$.
For $j=k$ this is the posterior $\bspi(X_k\mid Y^{[k]})$, since $Y^{(k,k]}=\bsone$.

\begin{theorem}\label{thm: PS}
	Let the notations be as above.
	\begin{enumerate}[label=(\roman*)]
		\item (Prediction) For $1\leq k\leq j\leq n$, we have
		\[
			\bspi(X_j\mid Y^{[k]})=\pr_{X_j}\circ\bspi(Y^{(k,j]}\times X_j\mid Y^{[k]}),
		\]
		and if $j<n$, we have
		\[
			\bspi(Y^{(k,j+1]}\times X_{j+1}\mid Y^{[k]})
			=(\pr_{Y^{(k,j]}},b_j)\circ\bspi(Y^{(k,j]}\times X_j\mid Y^{[k]})^\wedge
		\]
		for $\bspi_{Y^{[k]}}(s)$-almost every $(y_1,\ldots,y_k)\in Y^{[k]}$ for any $s\in S$.
		Here $b_j$ is regarded as a stochastic morphism on $S\times Y^{[k]}\times Y^{(k,j]}\times X_j=S\times Y^{[j]}\times X_j$.
		\item (Smoothing) For $1\leq j < k\leq n$, we have
		\[
			\bspi(X_j\mid Y^{[k]})=a_j^{\bspi(X_j\mid Y^{[j]})}\circ\bspi(X_{j+1}\mid Y^{[k]})^\wedge
		\]
		for $\bspi_{Y^{[k]}}(s)$-almost every $(y_1,\ldots,y_k)\in Y^{[k]}$ for any $s\in S$,
		where $a_j^{\bspi(X_j\mid Y^{[j]})}\colon S\times Y^{[j]}\times X_{j+1}\rightarrow X_j$
		is the Bayesian inverse of $a_j$ for the prior $\bspi(X_j\mid Y^{[j]})$,
		with $V=Y^{[j]}$ in \cref{def: BI V},
		regarded as a stochastic morphism on $S\times Y^{[k]}\times X_{j+1}$.
	\end{enumerate}
\end{theorem}

\begin{proof}
	The first assertion of \textit{(i)} is \cref{lem: properties}~$\ref{CL: P}$ for $a=\pr_{X_j}\colon Y^{(k,j]}\times X_j\rightarrow X_j$.
	For the second, we regard $\pi_j=\bspi_{Y^{[j]}\times X_j}$ as a joint law on $Y^{[k]}\times(Y^{(k,j]}\times X_j)$,
	and let $a\coloneqq(\pr_{Y^{(k,j]}},b_j)\colon S\times Y^{[j]}\times X_j\rightarrow Y^{(k,j+1]}\times X_{j+1}$.
	Since $Y^{[j]}=Y^{[k]}\times Y^{(k,j]}$, we have $(\pr_{Y^{[k]}},a)\circ\pi_j^\wedge=(\pr_{Y^{[j]}},b_j)\circ\pi_j=\pi_{j+1}\colon S\rightarrow Y^{[j+1]}\times X_{j+1}$
	by \eqref{eq: product} and \eqref{eq: recursion}.
	Hence the general case of \cref{lem: properties}~$\ref{CL: P}$, with $Y^{[k]}$ and $Y^{(k,j]}\times X_j$ in place of $X$ and $Y$, gives
	\[
		(\pr_{Y^{(k,j]}},b_j)\circ\bspi(Y^{(k,j]}\times X_j\mid Y^{[k]})^\wedge=\pi_{j+1}(Y^{(k,j+1]}\times X_{j+1}\mid Y^{[k]})
		=\bspi(Y^{(k,j+1]}\times X_{j+1}\mid Y^{[k]}),
	\]
	which proves the second assertion.
	The recursion of \textit{(i)} constructs $\bspi(Y^{(k,j]}\times X_j\mid Y^{[k]})$ from the posterior
	$\bspi(X_k\mid Y^{[k]})$, which in particular shows the existence of these conditional laws.

	We next prove \textit{(ii)}.  Let $1\leq j<k\leq n$.
	We first show that
	\begin{equation}\label{eq: PS1}
		a_j^{\bspi(X_j\mid Y^{[j]})}=\bspi(X_j\mid Y^{[j]}\times X_{j+1}).
	\end{equation}
	By construction of $\bspi$, the extension $\bspi_{Y^{[j]}\times X_{j+1}\times X_j}$ of $\pi_j$ is
	$
		\bspi_{Y^{[j]}\times X_{j+1}\times X_j}=(\pr_{Y^{[j]}},(a_j,\pr_{X_j}))\circ\pi_j,
	$
	hence by \cref{lem: properties}~$\ref{CL: P}$ applied to $\pi_j$ and $a=(a_j,\pr_{X_j})$,
	\[
		\kappa\coloneqq(a_j,\pr_{X_j})\circ\bspi(X_j\mid Y^{[j]})^\wedge
		=\bspi(X_{j+1}\times X_j\mid Y^{[j]}).
	\]
	By \cref{prop: construct V} with $V=Y^{[j]}$, the prior $\bspi(X_j\mid Y^{[j]})$
	and $a_j\colon Y^{[j]}\times X_j\rightarrow X_{j+1}$, we have
	$a_j^{\bspi(X_j\mid Y^{[j]})}=\kappa(X_j\mid X_{j+1})$,
	and since $\kappa=\bspi(X_{j+1}\times X_j\mid Y^{[j]})$, the tower property \cref{thm: tower} gives
	$\bspi(X_j\mid Y^{[j]}\times X_{j+1})=\kappa(X_j\mid X_{j+1})$,
	which proves \eqref{eq: PS1}.

	We next show that
	\begin{equation}\label{eq: PS2}
		\bspi(X_j\mid Y^{[k]}\times X_{j+1})=\bspi(X_j\mid Y^{[j]}\times X_{j+1})\colon S\times Y^{[k]}\times X_{j+1}\rightarrow X_j.
	\end{equation}
	Let $\wt{Y}\coloneqq Y^{(j,k]}$, and let $b\coloneqq\pr_{\wt{Y}}\circ c_k\colon Y^{[j]}\times X_{j+1}\rightarrow \wt{Y}$,
	where $c_l\colon Y^{[j]}\times X_{j+1}\rightarrow Y^{(j,l]}\times X_l$ for $l=j+1,\ldots,k$
	is defined iteratively by $c_{j+1}\coloneqq(h_{j+1},\pr_{X_{j+1}})$ and
	$c_{l+1}\coloneqq(\pr_{Y^{(j,l]}},b_l)\circ c_l^\wedge$,
	with $b_l$ regarded as a stochastic morphism on $Y^{[j]}\times X_{j+1}\times Y^{(j,l]}\times X_l$.
	By construction of $\bspi$, the future observations $\wt{Y}=Y^{(j,k]}$ are generated from
	$Y^{[j]}\times X_{j+1}$ through $b$, independently of $X_j$, so that the extension
	$\bspi_{Y^{[j]}\times X_{j+1}\times \wt{Y}\times X_j}$ of $\bspi_{Y^{[j]}\times X_{j+1}\times X_j}$ is
	\[
		\bspi_{Y^{[j]}\times X_{j+1}\times \wt{Y}\times X_j}
		=(\pr_{Y^{[j]}\times X_{j+1}},\,b,\,\pr_{X_j})
		\circ\bspi_{Y^{[j]}\times X_{j+1}\times X_j}.
	\]
	Hence $X_j\perp_{\bspi}\wt{Y}\mid Y^{[j]}\times X_{j+1}$ by \cref{cor: CI}
	with $W=Y^{[j]}\times X_{j+1}$, $X=X_j$ and $Y=\wt{Y}$, and \eqref{eq: PS2} follows from \cref{lem: CI},
	noting that $Y^{[j]}\times X_{j+1}\times \wt{Y}=Y^{[k]}\times X_{j+1}$ up to the order of the factors.
	Finally, by the chain rule \cref{lem: properties}~$\ref{CL: C}$, \eqref{eq: PS2} and \eqref{eq: PS1}, we have
	\begin{align*}
		\bspi(X_j\mid Y^{[k]})
		&=\bspi(X_j\mid Y^{[k]}\times X_{j+1})\circ\bspi(X_{j+1}\mid Y^{[k]})^\wedge\\
		&=\bspi(X_j\mid Y^{[j]}\times X_{j+1})\circ\bspi(X_{j+1}\mid Y^{[k]})^\wedge
		=a_j^{\bspi(X_j\mid Y^{[j]})}\circ\bspi(X_{j+1}\mid Y^{[k]})^\wedge,
	\end{align*}
	which proves \textit{(ii)}.
\end{proof}

\begin{remark}
	Iterating \cref{thm: PS}~\textit{(i)} computes
	$\bspi(Y^{(k,j]}\times X_j\mid Y^{[k]})$ for $k<j\leq n$ from the filter $\bspi(X_k\mid Y^{[k]})$,
	by propagating the hidden state together with the future observations, and the prediction
	$\bspi(X_j\mid Y^{[k]})$ is its marginal.
	In the case without feedback, we have
	$\bspi(X_{j+1}\mid Y^{[k]})=a_j\circ \bspi(X_j\mid Y^{[k]})$ by \textit{(i)}, hence
	\[
		\bspi(X_j\mid Y^{[k]})=a_{j-1}\circ \cdots \circ a_k\circ \bspi(X_k\mid Y^{[k]}),\qquad k<j\leq n .
	\]
	Similarly, iterating \cref{thm: PS}~\textit{(ii)}, the smoothed posteriors $\bspi(X_j\mid Y^{[k]})$ for $j<k$ are
	computed backward from the filter by composing with the Bayesian inverses of the transitions.
	By \eqref{eq: PS1}, these are the conditional laws
	$\bspi(X_j\mid Y^{[j]}\times X_{j+1})$,
	the transition kernels of the hidden state process run backward, and
	\eqref{eq: PS2} says that the future observations carry no further information on $X_j$ once $X_{j+1}$ is known.
\end{remark}

\begin{remark}
	The whole section uses only the three rules $\ref{aug: P}$--$\ref{aug: I}$ of the augmentation
	given in \cref{lem: augmentations},
	through \cref{lem: properties,lem: CI,cor: CI} and the tower property \cref{thm: tower}.
	In particular, no densities and no explicit integrals appear. Measure theory is confined to
	the existence and uniqueness of the conditional law and to the verification of the three rules.
\end{remark}

%
%
%
\section{The Affine Gaussian Case: Kalman Filter and Smoother}\label{sec: Gaussian}
%
%
%

In this section, we apply the results of \S\ref{sec: filter} to the affine Gaussian model,
and recover the Kalman filter \cite{Kal60} and the Rauch--Tung--Striebel smoother \cite{RTS65},
for models with feedback, whose coefficients may depend on the past observations.
The key observation is that the conditional laws of a Gaussian joint law are again Gaussian,
with mean and covariance given by explicit matrix formulas.
Hence the identities of \cref{thm: BF,thm: PS}, which hold for conditional laws, become recursions for means and covariances.
All calculations reduce to the independence of the noise, \cref{lem: sum}, and two properties of the Gaussian family, \cref{lem: Gaussian family}, through the composition rule \cref{lem: composition}.

We say that a measurable space $X$ is \emph{Euclidean} if $X=\R^{d_X}$ for an integer $d_X\geq0$, called its dimension,
with its standard inner product and Borel $\sigma$-algebra, where $\R^0=\bsone$.
We denote by $M_{p\times q}(\R)$ the space of real $p\times q$ matrices and by $\Sym^+_d(\R)$ the space of
symmetric positive semidefinite $d\times d$ matrices,
and we also equip $M_{p\times q}(\R)$ and $\Sym^+_d(\R)$
with their Borel $\sigma$-algebras.
For $P\in\Sym^+_d(\R)$, we denote by $P^+$ the Moore--Penrose inverse of $P$, so that $P^+\in\Sym^+_d(\R)$,
$PP^+P=P$ and $P^+PP^+=P^+$.
Since $P^+=\lim_{\ve\downarrow0}(P^2+\ve I)^{-1}P$ for $P\in\Sym^+_d(\R)$, the map $P\mapsto P^+$ is measurable.

%
\subsection{Mean and Covariance}\label{sec: mean}
%

Let $\pi\colon S\rightarrow X$ be a stochastic morphism, and let $\bsg\colon X\rightarrow\R^{d'}$ be a measurable function
which is $\pi(s)$-integrable for any $s\in S$.
Note that we use boldface to emphasize higher dimensional Euclidean vectors.
We define the \emph{expectation} of $\bsg$ with respect to $\pi$ to be the measurable function
\[
	\bbE_\pi[\bsg]\colon S\rightarrow\R^{d'},\qquad
	\bbE_\pi[\bsg](s)\coloneqq\int_X\bsg(\bsx)\,\pi(s)(d\bsx).
\]
For $\bsg$ bounded, $\bbE_\pi[\bsg]=\pi^*\bsg$ is the pull-back of \cref{def: pullback} with $S$ in place of $X$.

\begin{remark}\label{rem: X}
	In formulas for functions and stochastic morphisms on a product of measurable spaces
	with a Euclidean factor $X$, we write $X$ for the coordinate projection $\pr_X$ to $X$,
	and for $\id_X$ when the product is $X$ itself, in accordance with the classical notation
	$\bbE[X]$ for a random variable $X$ (see \cref{rem: prX}, as well as 
	the paragraph at the end of Appendix \ref{sec: appendix B}).
\end{remark}

We define the \emph{mean} and the \emph{covariance matrix} of $\pi$ by 
$
	\bbE_\pi[X]\coloneqq\bbE_\pi[\id_X]\colon S\rightarrow X
$ 
and
\[
	\bbV_\pi[X](s)\coloneqq\int_X\bigl(\bsx-\bbE_\pi[X](s)\bigr)\bigl(\bsx-\bbE_\pi[X](s)\bigr)^{\top}\pi(s)(d\bsx),
\]
whenever $\bbE_\pi[\abs{X}^2]\coloneqq\bbE_\pi[\abs{\id_X}^2]$ is finite for every $s\in S$.  Then $\bbV_\pi[X]\colon S\rightarrow\Sym^+_{d_X}(\R)$ is measurable.
For a joint law $\pi\colon S\rightarrow X\times Y$ on Euclidean $X$ and $Y$, we have
$\bbE_\pi[X]=\bbE_{\pi_X}[X]$, and define the \emph{cross-covariance}
$\Cov_\pi[X,Y]\colon S\rightarrow M_{d_X\times d_Y}(\R)$
of $X$ and $Y$ by
\[
	\Cov_\pi[X,Y](s)\coloneqq\int_{X\times Y}\bigl(\bsx-\bbE_\pi[X](s)\bigr)\bigl(\bsy-\bbE_\pi[Y](s)\bigr)^{\top}\pi(s)(d\bsx,d\bsy).
\]
The mean and covariance matrix of $\pi$ are then the block matrices
\begin{equation}\label{eq: block}
	\bbE_\pi[X\times Y]=\begin{pmatrix}\bbE_{\pi_X}[X]\\ \bbE_{\pi_Y}[Y]\end{pmatrix},\qquad
	\bbV_\pi[X\times Y]=\begin{pmatrix}
		\bbV_{\pi_X}[X] & \Cov_\pi[X,Y]\\
		\Cov_\pi[X,Y]^{\top} & \bbV_{\pi_Y}[Y]
	\end{pmatrix}.
\end{equation}

%
\subsection{Gaussian Stochastic Morphisms}
%

For $\bsm\in\R^d$ and $P\in\Sym^+_d(\R)$, we denote by $\sN(\bsm,P)$ the Gaussian probability measure on $\R^d$
with mean $\bsm$ and covariance matrix $P$, characterized by its characteristic function
\begin{equation}\label{eq: char}
	\int_{\R^d}e^{i\bsu^{\top}\bsx}\sN(\bsm,P)(d\bsx)=\exp\bigl(i\bsu^{\top}\bsm-\tfrac12\bsu^{\top}P\bsu\bigr),\qquad \bsu\in\R^d.
\end{equation}
We allow the covariance matrix $P$ to be singular.
Let
\[
	\Theta_d\coloneqq\R^d\times\Sym^+_d(\R),
\]
which is a commutative monoid under addition.
Since $\sN(\bsm,P)$ is the image of $\sN(0,I)$ under the affine map $\bsx\mapsto\bsm+P^{1/2}\bsx$,
and $P^{1/2}$ depends continuously on $P$, the Gaussian family
\[
	\sN\colon\Theta_d\rightarrow\R^d,\qquad(\bsm,P)\mapsto\sN(\bsm,P)
\]
is a stochastic morphism, and for a measurable space $S$ and measurable functions
$\bsm\colon S\rightarrow\R^d$ and $P\colon S\rightarrow\Sym^+_d(\R)$,
the composition $\sN(\bsm,P)\coloneqq\sN\circ(\bsm,P)\colon S\rightarrow\R^d$, $s\mapsto\sN(\bsm(s),P(s))$,
is a stochastic morphism.

We next introduce the sum and the linear image of stochastic morphisms with values in a Euclidean space,
which express adding an independent noise and applying a linear map.
Let $S$ be a measurable space.
For stochastic morphisms $p,q\colon S\rightarrow\R^d$, we define the \emph{sum}
\[
	p+q\coloneqq{+}\circ(p,q)\colon S\rightarrow\R^d,
\]
where $(p,q)\colon S\rightarrow\R^d\times\R^d$ is the product of \S\ref{sec: CL}
and ${+}\colon\R^d\times\R^d\rightarrow\R^d$ is the stochastic morphism induced by the measurable
addition map $(\bsx,\bsx')\mapsto\bsx+\bsx'$.
Since $(p,q)$ is the product, $p+q$ is the law of the sum of independent variables.
In particular, $p+p\neq2p$ in general.

\Cref{lem: sum} below states that a noise which depends only on the parameter is independent of the input,
which is the property of the noise used in Kalman filtering.

\begin{lemma}[Independent Noise]\label{lem: sum}
	Let $\kappa\colon S\rightarrow Y$ be a stochastic morphism between measurable spaces.
	Furthermore, let $\bsf\colon S\times Y\rightarrow\R^d$ be a measurable map, and let
	$\varepsilon\colon S\rightarrow\R^d$ be a stochastic morphism,
	which we view as a morphism on $S\times Y$.
	Then
	\[
		(\bsf+\varepsilon)\circ\kappa^\wedge=\bsf\circ\kappa^\wedge+\varepsilon\colon S\rightarrow\R^d.
	\]
\end{lemma}

\begin{proof}
	We regard $\varepsilon$ as the morphism $\varepsilon\circ\pr_S\colon S\times Y\rightarrow\R^d$.
	By \eqref{eq: product}, we have 
	\[
		(\bsf,\varepsilon\circ\pr_S)=(\bsf\times\id_{\R^d})\circ(\varepsilon\circ\pr_S)^\wedge
		\colon S\times Y\rightarrow\R^d\times\R^d,
	\] 
	where 
	$(\varepsilon\circ\pr_S)^\wedge=(\id_{S\times Y},\varepsilon\circ\pr_S)\colon S\times Y\rightarrow S\times Y\times\R^d$ is the augmentation with parameter space 
	$S\times Y$. By \cref{lem: augmentations}~$\ref{aug: I}$ applied to $\kappa$ and $\varepsilon$, we have $(\varepsilon\circ\pr_S)^\wedge\circ\kappa^\wedge=(\id_S,\kappa,\varepsilon)=(\kappa^\wedge,\varepsilon)$. Hence by \eqref{eq: product},
	\[
		(\bsf,\varepsilon\circ\pr_S)\circ\kappa^\wedge=(\bsf\times\id_{\R^d})\circ(\kappa^\wedge,\varepsilon)=(\bsf\circ\kappa^\wedge,\;\varepsilon),
	\]
	and composing with ${+}$ gives the assertion.
\end{proof}

For a measurable function $A\colon S\rightarrow M_{d'\times d}(\R)$,
we define the \emph{linear image} of a stochastic morphism
$p\colon S\rightarrow\R^d$ to be the stochastic morphism
\[
	Ap\colon S\rightarrow\R^{d'},\qquad (Ap)(s)\coloneqq A(s)(p(s)),
\]
where $A(s)(p(s))$ is the push-forward of the probability measure $p(s)$ on $\R^d$
under the linear map $A(s)\colon\R^d\rightarrow\R^{d'}$.
The Gaussian family has the following two properties.

\begin{lemma}[Gaussian Family]\label{lem: Gaussian family}
	Let $S$ be a measurable space, and let $\bsm,\bsm'\colon S\rightarrow\R^d$,
	$P,P'\colon S\rightarrow\Sym^+_d(\R)$
	and $A\colon S\rightarrow M_{d'\times d}(\R)$ be measurable.
	\begin{enumerate}[label=(\roman*)]
		\item\label{GF: conv} (Sum) $\sN(\bsm,P)+\sN(\bsm',P')=\sN(\bsm+\bsm',\,P+P')$.
		\item\label{GF: equiv} (Linear image) $A\,\sN(\bsm,P)=\sN(A\bsm,\,APA^{\top})$.
	\end{enumerate}
	In particular, we have $\sN(\bsm,P)=\bsm+\sN(0,P)$.
\end{lemma}

\begin{proof}
	It suffices to prove the assertions for each $s\in S$.
	The characteristic function of the image of a product measure under the addition is the product of the characteristic functions,
	and that of $A(s)\bigl(\sN(\bsm(s),P(s))\bigr)$ at $\bsu$ is that of $\sN(\bsm(s),P(s))$ at $A(s)^{\top}\bsu$.
	Both assertions follow from \eqref{eq: char} and the uniqueness of characteristic functions.
	The last statement is $\ref{GF: conv}$ applied to $\bsm=\sN(\bsm,0)$ and $\sN(0,P)$.
\end{proof}

In the language of stochastic morphisms, $\ref{GF: conv}$ states that
we have an equality of stochastic morphisms
${+}\circ(\sN\times\sN)=\sN\circ{+}\colon\Theta_d\times\Theta_d\rightarrow\R^d$,
that is, the Gaussian family $\sN\colon\Theta_d\rightarrow\R^d$ is a homomorphism of monoids in the stochastic category,
and $\ref{GF: equiv}$ states that $\sN$ is compatible with linear maps.

\begin{definition}\label{def: Gaussian}
	Let $S$ be a measurable space, and let $X$ and $Y$ be Euclidean.
	\begin{enumerate}[label=(\roman*)]
		\item We say that a stochastic morphism $\mu\colon S\rightarrow Y$
		is \emph{Gaussian} if $\mu(s)$ is a Gaussian measure for any $s\in S$, that is,
		\[
			\mu=\sN(\bsm,P)\colon S\rightarrow Y
		\]
		for measurable functions $\bsm\colon S\rightarrow Y$ and $P\colon S\rightarrow\Sym^+_{d_Y}(\R)$,
		which are the mean $\bbE_\mu[Y]$ and the covariance matrix $\bbV_\mu[Y]$ of $\mu$.
		\item We say that a stochastic morphism $h\colon S\times X\rightarrow Y$ is
		\emph{affine Gaussian in $X$} if
		\[
			h=\sN(AX+\bsb,\,R)\colon S\times X\rightarrow Y,
		\]
		that is, 
		there exist measurable functions
		$A\colon S\rightarrow M_{d_Y\times d_X}(\R)$, $\bsb\colon S\rightarrow Y$, $R\colon S\rightarrow\Sym^+_{d_Y}(\R)$ such that
		$h(s,\bsx)=\sN(A(s)\bsx+\bsb(s),\,R(s))$ 
		for any $(s,\bsx)\in S\times X$.			
		Here, $A$, $\bsb$ and $R$ are regarded 
		as functions on $S\times X$ through $\pr_S$.
		We call $AX+\bsb\colon S\times X\rightarrow Y$ the \emph{mean function}, 
		$A$ the \emph{linear part},
		and $R$ the \emph{covariance matrix} of $h$.
		Thus $h$ may depend on the parameter $s\in S$ in an arbitrary measurable way,
		while it is affine in $X$ with covariance matrix independent of $X$.
		Note that here, $X$ denotes the coordinate projection $\pr_X\colon S\times X\rightarrow X$
		as in \cref{rem: X}.
		\end{enumerate}
\end{definition}

The following composition rule is a consequence of \cref{lem: sum,lem: Gaussian family}.

\begin{proposition}[Composition of Affine Gaussian Morphisms]\label{lem: composition}
	Let $S$ be a measurable space, and let $X$, $Y$ and $Z$ be Euclidean.
	\begin{enumerate}[label=(\roman*)]
		\item\label{comp: rule} (Composition)
		For stochastic morphisms
		\begin{align*}
			h_1&=\sN(A_1X+\bsb_1,\,R_1)\colon S\times X\rightarrow Y,&
				h_2&=\sN(A_2Y+\bsb_2,\,R_2)\colon S\times Y\rightarrow Z,
		\end{align*}
		affine Gaussian in $X$ and in $Y$ respectively,
		the composition $h_2\circ h_1^\wedge\colon S\times X\rightarrow Z$ is affine Gaussian in $X$, given by
		\begin{equation}\label{eq: composition rule}
			h_2\circ h_1^\wedge=\sN\bigl(A_2A_1X+A_2\bsb_1+\bsb_2,\;A_2R_1A_2^{\top}+R_2\bigr).
		\end{equation}
		\item\label{comp: push} (Push-forward) For a Gaussian $\mu=\sN(\bsm,P)\colon S\rightarrow X$ and
		$h=\sN(AX+\bsb,\,R)\colon S\times X\rightarrow Y$ affine Gaussian in $X$, we have
		\begin{equation}\label{eq: pushforward rule}
			h\circ\mu^\wedge=\sN\bigl(A\bsm+\bsb,\;APA^{\top}+R\bigr)\colon S\rightarrow Y.
		\end{equation}
		\item\label{comp: joint} (Joint law) For $\mu$ and $h$ as in \ref{comp: push}, the stochastic morphism
		$(\pr_X,h)\colon S\times X\rightarrow X\times Y$ is affine Gaussian in $X$, and the joint law
		$\kappa\coloneqq(\pr_X,h)\circ\mu^\wedge\colon S\rightarrow X\times Y$ is Gaussian with
		\begin{equation}\label{eq: joint law}
			\bbE_\kappa[X\times Y]=\begin{pmatrix}\bsm\\ A\bsm+\bsb\end{pmatrix},\qquad
			\bbV_\kappa[X\times Y]=
			\begin{pmatrix}
				P & PA^{\top} \\
				AP &APA^{\top}+R
			\end{pmatrix}.
		\end{equation}
	\end{enumerate}
\end{proposition}

\begin{proof}
	For $\ref{comp: rule}$, we have $h_2=A_2Y+\bsb_2+\sN(0,R_2)$ by \cref{lem: Gaussian family}, so that
	regarding $h_2$ as a stochastic morphism on $S\times X\times Y$, we have
	$h_2=\bsf+\sN(0,R_2)\circ\pr_{S\times X}\colon S\times X\times Y\rightarrow Z$ for the measurable function
	$\bsf\colon S\times X\times Y\rightarrow Z$, $\bsf(s,\bsx,\bsy)\coloneqq A_2(s)\bsy+\bsb_2(s)$.
	Hence by \cref{lem: sum} applied with $S\times X$ in place of $S$, we have
	\[
		h_2\circ h_1^\wedge=\bsf\circ h_1^\wedge+\sN(0,R_2)=A_2h_1+\bsb_2+\sN(0,R_2).
	\]
	Since $h_1=A_1X+\bsb_1+\sN(0,R_1)$ by \cref{lem: Gaussian family},
	\cref{lem: Gaussian family}~$\ref{GF: equiv}$ gives
	\[
		A_2h_1=A_2A_1X+A_2\bsb_1+A_2\sN(0,R_1)=A_2A_1X+A_2\bsb_1+\sN(0,A_2R_1A_2^{\top}),
	\]
	and \eqref{eq: composition rule} follows from \cref{lem: Gaussian family}~$\ref{GF: conv}$.
	Furthermore,
	$\ref{comp: push}$ is $\ref{comp: rule}$ with $X$, $Y$ and $Z$ replaced by $\bsone$, $X$ and $Y$,
	applied to $h_1=\mu$ regarded as $S\times\bsone\rightarrow X$ with $A_1=0$, $\bsb_1=\bsm$ and $R_1=P$, and $h_2=h$.
	Finally,
	for $\ref{comp: joint}$, 
	let $\iota_Y\colon Y\rightarrow X\times Y$ be the linear map $\bsy\mapsto(0,\bsy)$.
	Since $(\pr_X,h)(s,\bsx)$ is the image of $h(s,\bsx)$ under the map $\bsy\mapsto(\bsx,0)+\iota_Y(\bsy)$,
	\cref{lem: Gaussian family}~$\ref{GF: equiv}$ gives
	\[
		(\pr_X,h)=\sN\biggl(\begin{pmatrix}I\\A\end{pmatrix}X+\begin{pmatrix}0\\ \bsb\end{pmatrix},\;\begin{pmatrix}0&0\\0&R\end{pmatrix}\biggr),
	\]
	which is affine Gaussian in $X$.  Then $\ref{comp: push}$ gives \eqref{eq: joint law}.
\end{proof}

\begin{remark}
	For each measurable space $S$, the stochastic morphisms $S\times X\to Y$ affine Gaussian in $X$, for Euclidean $X$ and $Y$, are closed under the composition of \cref{lem: composition}~$\ref{comp: rule}$. For $S=\bsone$, these are the morphisms of the Markov category \textbf{Gauss} of \cite[\S6]{F20}, and the computations of this section take place in this parameterized version of \textbf{Gauss}.
	The inclusion is strict: a linear part or covariance depending on $S$, or a non-Euclidean $S$, gives morphisms which are not in \textbf{Gauss}, and these are what the Kalman filter with feedback requires.
\end{remark}

We next show that a Gaussian joint law has a Gaussian conditional law, given by an explicit formula.
The linear algebra involved is the following.

\begin{lemma}\label{lem: Schur}
	Consider the matrix
	\[
		M=\begin{pmatrix}P_X & C\\ C^{\top} & P_Y\end{pmatrix}\in\Sym^+_{d_X+d_Y}(\R)
	\]
	with $P_X\in\Sym^+_{d_X}(\R)$, $P_Y\in\Sym^+_{d_Y}(\R)$ and $C\in M_{d_X\times d_Y}(\R)$, and let
	\[
		K\coloneqq CP_Y^+\in M_{d_X\times d_Y}(\R),\qquad \Sigma\coloneqq P_X-CP_Y^+C^{\top}.
	\]
	Then $KP_Y=C$, and $\Sigma=P_X-KP_YK^{\top}=(I\;\,{-K})\,M\,(I\;\,{-K})^{\top}\in\Sym^+_{d_X}(\R)$.
\end{lemma}

\begin{proof}
	If $P_Y\bsw=0$, then $0\leq(\bsv,t\bsw)M(\bsv,t\bsw)^{\top}=\bsv^{\top}P_X\bsv+2t\,\bsv^{\top}C\bsw$
	for all $\bsv\in\R^{d_X}$ and $t\in\R$, hence $C\bsw=0$.
	Thus $\Ker P_Y\subset\Ker C$, so $\Im C^{\top}\subset\Im P_Y$, and since $P_Y^+P_Y=P_YP_Y^+$ is the orthogonal
	projection onto $\Im P_Y$, we have $KP_Y=CP_Y^+P_Y=C$.
	Then $(I\;\,{-K})M(I\;\,{-K})^{\top}=P_X-KC^{\top}-CK^{\top}+KP_YK^{\top}=P_X-CP_Y^+C^{\top}=\Sigma$,
	since $KC^{\top}=CP_Y^+C^{\top}=CK^{\top}$ and $KP_YK^{\top}=CK^{\top}$.
	In particular $\Sigma=P_X-KP_YK^{\top}$ is positive semidefinite.
\end{proof}

\begin{proposition}[Gaussian Conditioning]\label{prop: Gaussian conditioning}
	Let
	\[
		\pi\colon S\rightarrow X\times Y
	\]
	be a Gaussian joint law on Euclidean $X$ and $Y$ such that
	\begin{align*}
		\bbE_\pi[X\times Y]&=
		\begin{pmatrix} \bsm_X\\ \bsm_Y
		\end{pmatrix},  &
		\bbV_\pi[X\times Y]&=
		\begin{pmatrix}
			P_X  & C\\
			C^\top & P_Y
		\end{pmatrix},
	\end{align*}
	where $C\coloneqq\Cov_\pi[X,Y]\colon S\rightarrow M_{d_X\times d_Y}(\R)$ is the cross-covariance,
	and let $K\coloneqq CP_Y^+\colon S\rightarrow M_{d_X\times d_Y}(\R)$.
	Then the conditional law $\pi(X\mid Y)\colon S\times Y\rightarrow X$ of $X$ given $Y$ for $\pi$ satisfies
	\[
		\pi(X\mid Y)=\sN\bigl(\bsm_X+K(Y-\bsm_Y),\;P_X-KP_YK^{\top}\bigr)
	\]
	for $\pi_Y(s)$-almost every $\bsy\in Y$ for any $s\in S$,
	where the mean and covariance are evaluated at $s$.
\end{proposition}

\begin{proof}
	Let $\Sigma\coloneqq P_X-KP_YK^{\top}$.  The functions $K$ and $\Sigma$ on $S$ are measurable,
	since the moments and the Moore--Penrose inverse are measurable.
	By \cref{lem: Schur} applied to the covariance matrix of $\pi(s)$ for each $s\in S$,
	$\Sigma$ takes values in $\Sym^+_{d_X}(\R)$, hence
	\[
		\kappa\coloneqq\sN\bigl(\bsm_X+K(Y-\bsm_Y),\,\Sigma\bigr)\colon S\times Y\rightarrow X
	\]
	is a stochastic morphism affine Gaussian in $Y$ with linear part $K$, which is the right hand side of the assertion.
	By the uniqueness of the conditional law, \cref{thm: CL}, it suffices to show that $\kappa$ is a conditional law
	of $X$ given $Y$ for $\pi$, that is, $(\pr_Y,\kappa)\circ(\pi_Y)^\wedge=\pi_{Y\times X}$.
	By \cref{lem: composition}~$\ref{comp: joint}$ for the Gaussian $\pi_Y=\sN(\bsm_Y,P_Y)\colon S\rightarrow Y$
	and $\kappa$, which is affine Gaussian in $Y$, the left hand side is Gaussian with mean
	$(\bsm_Y,\,\bsm_X+K(\bsm_Y-\bsm_Y))=(\bsm_Y,\bsm_X)$ and covariance matrix
	\[
		\begin{pmatrix}P_Y & P_YK^{\top}\\ KP_Y & KP_YK^{\top}+\Sigma\end{pmatrix}
		=\begin{pmatrix}P_Y & C^{\top}\\ C & P_X\end{pmatrix},
	\]
	by $KP_Y=C$ of \cref{lem: Schur}.
	These are the mean and covariance matrix of the Gaussian measure $\pi_{Y\times X}(s)$,
	hence the two sides agree.
\end{proof}

Combining \cref{lem: composition,prop: Gaussian conditioning} with \cref{prop: construct V}, we see that
the Bayesian inverse of an observation affine Gaussian in the hidden state, with respect to a Gaussian prior,
is affine Gaussian in the observation,
with coefficients given by the Kalman gain and the posterior covariance.
In other words, affine Gaussian stochastic morphisms are closed under Bayesian inversion.

\begin{corollary}[Affine Gaussian Bayesian Inverse]\label{cor: Gaussian BI}
	Let $V$ be a measurable space, and let $X$ and $Y$ be Euclidean.
	Consider a Gaussian prior and an observation affine Gaussian in $X$
	\[
		\mu=\sN(\bsm,P)\colon S\times V\rightarrow X,\qquad
		h=\sN(AX+\bsb,\,R)\colon V\times X\rightarrow Y,
	\]
	as in \cref{def: BI V}, where $\bsm,P$ are measurable functions on $S\times V$
	and $A,\bsb,R$ are measurable functions on $V$, regarded as functions on $S\times V$ through $\pr_V$.
	Let
	\[
		D\coloneqq APA^{\top}+R\colon S\times V\rightarrow\Sym^+_{d_Y}(\R),\qquad
		K\coloneqq PA^{\top}D^+\colon S\times V\rightarrow M_{d_X\times d_Y}(\R).
	\]
	Then the law of the observation and the Bayesian inverse of $h$ with respect to $\mu$ are
	\begin{align*}
		h\circ\mu^\wedge&=\sN(A\bsm+\bsb,\,D)\colon S\times V\rightarrow Y,\\
		h^\mu&=\sN\bigl(\bsm+K(Y-A\bsm-\bsb),\;P-KDK^{\top}\bigr)\colon S\times V\times Y\rightarrow X,
	\end{align*}
	the second for $(h\circ\mu^\wedge)(s,v)$-almost every $\bsy\in Y$ for any $(s,v)\in S\times V$.
\end{corollary}

\begin{proof}
	The first assertion is \eqref{eq: pushforward rule} with parameter space $S\times V$.
	Let $\kappa\coloneqq(\pr_X,h)\circ\mu^\wedge\colon S\times V\rightarrow X\times Y$.
	By \eqref{eq: product}, $\kappa_{Y\times X}=(h,\pr_X)\circ\mu^\wedge$, hence $h^\mu=\kappa(X\mid Y)$ by \cref{prop: construct V}.
	By \eqref{eq: joint law}, $\kappa$ is Gaussian with $\bsm_X=\bsm$, $\bsm_Y=A\bsm+\bsb$, $P_X=P$, $P_Y=D$ and
	$C=PA^{\top}$, so that $CP_Y^+=K$, and the second assertion follows from \cref{prop: Gaussian conditioning}.
\end{proof}

%
\subsection{Kalman Filter and Rauch--Tung--Striebel Smoother}
%

We now consider the model of \cref{def: HMP} in the affine Gaussian case.
Let $X_i$ and $Y_i$ for $1\leq i\leq n$ be Euclidean of dimensions $d_X$ and $d_Y$ independent of $i$,
with Gaussian prior
\[
	\mu=\sN(\bsm_{1\mid0},P_{1\mid0})\colon S\rightarrow X_1,
\]
and transitions and observations affine Gaussian in $X_i$
\begin{align*}
	a_i&=\sN(A_iX_i+\bsb_i,\,Q_i)\colon Y^{[i]}\times X_i\rightarrow X_{i+1},& 1\leq i<n,\\
	h_i&=\sN(H_iX_i+\bsc_i,\,R_i)\colon Y^{[i-1]}\times X_i\rightarrow Y_i,& 1\leq i\leq n,
\end{align*}
where the coefficients
\begin{align*}
	A_i&\colon Y^{[i]}\rightarrow M_{d_X\times d_X}(\R),& \bsb_i&\colon Y^{[i]}\rightarrow X_{i+1},& Q_i&\colon Y^{[i]}\rightarrow\Sym^+_{d_X}(\R),\\
	H_i&\colon Y^{[i-1]}\rightarrow M_{d_Y\times d_X}(\R),& \bsc_i&\colon Y^{[i-1]}\rightarrow Y_i,& R_i&\colon Y^{[i-1]}\rightarrow\Sym^+_{d_Y}(\R)
\end{align*}
are measurable.
The coefficients may depend on the past observations, so that the model has feedback,
and the covariance matrices may be singular.
The classical Kalman filter is the case of coefficients not depending on the observations,
in which the model and its filter lie in the Markov category \textbf{Gauss} of \cite[\S6]{F20};
with feedback, the filter and the smoother are not morphisms of \textbf{Gauss},
since $\bsm_{k\mid k}$ and $P_{k\mid k}$ are no longer affine and constant in the observations.
Gaussian models whose coefficients depend on the past observations are known as conditionally Gaussian sequences,
and the corresponding filter and backward interpolation (smoothing) equations are classical,
see \cite[Theorems~13.4 and~13.12]{LS01}, where integrability of the coefficients is assumed
and the noises of the state and the observation may be correlated.

In our framework, \cref{thm: Kalman} is obtained from \cref{thm: BF,thm: PS} by the same argument
as in the case without feedback, for arbitrary measurable coefficients and uniformly in the parameter $s\in S$.
The means in \cref{thm: Kalman} are the means of the conditional laws, such as $\bspi(X_k\mid Y^{[k]})$,
which are defined without any integrability of $X_k$,
and coincide with the conditional expectations of $X_k$ given $Y^{[k]}$ when $X_k$ is integrable.
Let $\bspi\colon S\rightarrow Y^{[n]}\times X^{[n]}$ be the joint law of the model.
Functions on $Y^{[i]}$ are regarded as functions on $S\times Y^{[k]}$ for $i\leq k$ through the projection,
and we apply the results of this section with parameter space $S\times Y^{[k]}$ accordingly.

\begin{definition}[Kalman Recursion]\label{def: Kalman recursion}
	We define measurable functions $\bsm_{j\mid k}$ with values in $X_j$ and $P_{j\mid k}$ with values in
	$\Sym^+_{d_X}(\R)$ on $S\times Y^{[k]}$, for $0\leq k\leq n$ and $1\leq j\leq n$ with $j\leq k+1$, as follows.
	\begin{enumerate}[label=(\roman*)]
		\item (Update) For $1\leq k\leq n$, let
		$
			D_k\coloneqq H_kP_{k\mid k-1}H_k^{\top}+R_k
		$
		and
		$
			K_k\coloneqq P_{k\mid k-1}H_k^{\top}D_k^+
		$
		on $S\times Y^{[k-1]}$, and
		\[
			\bsm_{k\mid k}\coloneqq \bsm_{k\mid k-1}+K_k(Y_k-H_k\bsm_{k\mid k-1}-\bsc_k)
		\]
		and
		$
			P_{k\mid k}\coloneqq P_{k\mid k-1}-K_kD_kK_k^{\top}
		$
		on $S\times Y^{[k]}=S\times Y^{[k-1]}\times Y_k$.
		\item (Prediction) For $1\leq k<n$, let
		\[
			\bsm_{k+1\mid k}\coloneqq A_k\bsm_{k\mid k}+\bsb_k
		\]
		and
		$
			P_{k+1\mid k}\coloneqq A_kP_{k\mid k}A_k^{\top}+Q_k
		$
		on $S\times Y^{[k]}$.
		\item (Smoothing) For $1\leq j<k\leq n$, let
		$
			J_j\coloneqq P_{j\mid j}A_j^{\top}P_{j+1\mid j}^+
		$
		on $S\times Y^{[j]}$, and
		\[
			\bsm_{j\mid k}\coloneqq \bsm_{j\mid j}+J_j(\bsm_{j+1\mid k}-\bsm_{j+1\mid j})
		\]
		and
		$
			P_{j\mid k}\coloneqq P_{j\mid j}+J_j(P_{j+1\mid k}-P_{j+1\mid j})J_j^{\top}
		$
		on $S\times Y^{[k]}$.
	\end{enumerate}
	Here, $\bsm_{1\mid0}$ and $P_{1\mid0}$ are given by the prior, and the recursions run forward in $k$
	for (i)--(ii) and backward in $j$ for (iii).
\end{definition}

In the terminology of \cite{SS23}, the matrices $D_k$, $K_k$ and $J_j$ are the innovation covariance,
the Kalman gain and the smoother gain, and $D$ and $K$ of \cref{cor: Gaussian BI} are the same matrices for a single observation.

\begin{theorem}[Kalman Filter and Rauch--Tung--Striebel Smoother]\label{thm: Kalman}
	Let the notations be as in \cref{def: Kalman recursion}.
	\begin{enumerate}[label=(\roman*)]
		\item\label{Kalman: filter} (Filter) For $1\leq k\leq n$, we have
		$\bspi(X_k\mid Y^{[k]})=\sN(\bsm_{k\mid k},P_{k\mid k})$.
		\item\label{Kalman: predict} (Prediction) For $1\leq k<n$, the predicted prior satisfies
		$\mu_{k+1}=\bspi(X_{k+1}\mid Y^{[k]})=\sN(\bsm_{k+1\mid k},P_{k+1\mid k})$.
		\item\label{Kalman: smooth} (Smoothing) For $1\leq j<k\leq n$, we have
		$\bspi(X_j\mid Y^{[k]})=\sN(\bsm_{j\mid k},P_{j\mid k})$.
	\end{enumerate}
	The equalities hold for $\bspi_{Y^{[k]}}(s)$-almost every $(\bsy_1,\ldots,\bsy_k)\in Y^{[k]}$ for any $s\in S$.
\end{theorem}

\begin{proof}
	We choose the Gaussian versions of the conditional laws, which is allowed since the
	identities of \cref{thm: BF,thm: PS} hold for any choice.
	Let $\mu_1\coloneqq\mu$.  By \cref{prop: construct} for $k=1$ and \cref{thm: BF} for $k>1$, we have
	$\bspi(X_k\mid Y^{[k]})=h_k^{\mu_k}$, and by \cref{thm: PS}~\textit{(i)}, we have
	$\bspi(X_{k+1}\mid Y^{[k]})=a_k\circ\bspi(Y^{[k]}\times X_k\mid Y^{[k]})=\mu_{k+1}$.
	Hence $\ref{Kalman: filter}$ and $\ref{Kalman: predict}$ follow by induction on $k$.
	If $\mu_k=\sN(\bsm_{k\mid k-1},P_{k\mid k-1})$, then \cref{cor: Gaussian BI}
	with $V=Y^{[k-1]}$, $A=H_k$, $\bsb=\bsc_k$ and $R=R_k$ gives
	$h_k^{\mu_k}=\sN(\bsm_{k\mid k},P_{k\mid k})$,
	the innovation covariance $D_k$ and the Kalman gain $K_k$ of \cref{def: Kalman recursion}~(i)
	being $D$ and $K$ of \cref{cor: Gaussian BI}.
	Then \eqref{eq: pushforward rule} applied to
	$a_k$, regarded as a stochastic morphism on $S\times Y^{[k]}\times X_k$, with parameter space $S\times Y^{[k]}$ gives
	$\mu_{k+1}=\sN(\bsm_{k+1\mid k},P_{k+1\mid k})$.

	For $\ref{Kalman: smooth}$, we use backward induction on $j$, the case $j=k$ being 
	$\ref{Kalman: filter}$.
	By \cref{cor: Gaussian BI} for the prior $\bspi(X_j\mid Y^{[j]})=\sN(\bsm_{j\mid j},P_{j\mid j})$ and the
	observation $a_j$, whose law is $\sN(\bsm_{j+1\mid j},P_{j+1\mid j})$ by $\ref{Kalman: predict}$, we have
	\[
		a_j^{\bspi(X_j\mid Y^{[j]})}=\sN\bigl(\bsm_{j\mid j}+J_j(X_{j+1}-\bsm_{j+1\mid j}),\;P_{j\mid j}-J_jP_{j+1\mid j}J_j^{\top}\bigr)
	\]
	on $S\times Y^{[j]}\times X_{j+1}$.  By \cref{thm: PS}~\textit{(ii)}, $\bspi(X_j\mid Y^{[k]})$ is the
	composition of this stochastic morphism with $\bspi(X_{j+1}\mid Y^{[k]})=\sN(\bsm_{j+1\mid k},P_{j+1\mid k})$,
	which is $\sN(\bsm_{j\mid k},P_{j\mid k})$ by \eqref{eq: pushforward rule}.
\end{proof}

\begin{remark}
	The smoothing gain $J_j$ is the Kalman gain of \cref{cor: Gaussian BI} for the transition $a_j$ viewed as an observation
	of $X_j$ with prior $\bspi(X_j\mid Y^{[j]})$, in accordance with \eqref{eq: PS1},
	which identifies the Bayesian inverse of the transition with the backward transition
	$\bspi(X_j\mid Y^{[j]}\times X_{j+1})$ of the hidden state process.
	The smoother is thus the composition of the backward transition kernels with the filter, and the
	smoothing recursion of \cref{def: Kalman recursion} expresses this composition.
	In this way, \cref{thm: Kalman} is \cref{thm: BF,thm: PS} evaluated in the category of affine Gaussian
	stochastic morphisms.
\end{remark}

\begin{remark}[Extended Kalman Filter]
	Consider a nonlinear model with prior $\mu\colon\bsone\rightarrow X$,
	transitions $\sN(\bsf_k(X_k),Q_k)$ and observations
	$\sN(\bsg_k(X_k),R_k)$.  Then
	the extended Kalman filter linearizes $\bsf_k$ and $\bsg_k$ at the current filtered and predicted means,
	which are functions of the past observations.
	The resulting surrogate model is an instance of the model of this section with feedback,
	and \cref{thm: Kalman} applied to it gives the extended Kalman filter and smoother, exact for the surrogate model.
\end{remark}

\appendix
%
%
%
\section{Existence and Uniqueness of the Conditional Law}\label{sec: appendix A}
%
%
%

In this appendix, we prove \cref{thm: CL} concerning the existence and uniqueness of the conditional law.
This theorem is known.  The existence is proved in \cite[Theorem~1.25 (i)]{Kall17} when $X$ is standard Borel, and \cite[Theorem~2.29]{For21} for $X$ countably generated, with uniqueness as in \cite[Lemma~2.28]{For21}.   
See also \cite[Theorem~3.2]{Deg25}.
We give a complete proof for the convenience of the reader,
in our notation, and with the martingale argument written out.

For reference, we recall the definition of the conditional law and restate the theorem as \cref{thm: main}.
For a measurable space $X$, we denote by $\sB_X$ its $\sigma$-algebra.

\begin{theorem}\label{thm: main}
	Consider a stochastic morphism $\pi\colon S\rightarrow X\times Y$.
	 A \emph{conditional law} $\pi(Y\mid X)$ of $Y$ given $X$ for $\pi$ is a stochastic morphism
	\[
		\pi(Y\mid X)\colon S\times X\rightarrow Y,
	\]
	which satisfies
	\begin{equation}\label{eq: CL appendix}
		\pi=(\pr_X,\pi(Y\mid X))\circ(\id_S,\pi_X),
	\end{equation}
	where $\pi_X=\pr_X\circ\pi$ is the marginal.
	If we assume that $X$ is countably generated, then we have the following.
	\begin{enumerate}[label=(\roman*)]
		\item If $Y$ is countably generated, and if a conditional law
		 $\pi(Y\mid X)$ of $Y$ given $X$ for $\pi$ exists,
		then such $\pi(Y\mid X)$ is unique $\pi_X(s)$-almost everywhere for any $s\in S$.
		\item In particular, if $Y$ is standard Borel, then
		a conditional law
		 $\pi(Y\mid X)$ of $Y$ given $X$ for $\pi$ exists,
		 and is determined uniquely $\pi_X(s)$-almost everywhere for any $s\in S$.
	\end{enumerate}
\end{theorem}

Note that for $s\in S$, the measure $(\pi_X)^\wedge(s)$ on $S\times X$ is the product of the Dirac measure $\delta_s$
and $\pi_X(s)$. Hence for any stochastic morphism $\kappa\colon S\times X\rightarrow Y$, we have
\begin{equation}\label{eq: CL rectangle}
	((\pr_X,\kappa)\circ(\pi_X)^\wedge)(s)(A\times B)=\int_A \kappa(s,x)(B)\,\pi_X(s)(dx),\quad A\in\sB_X,\quad B\in\sB_Y.
\end{equation}
Since a probability measure on $X\times Y$ is determined by its values on measurable rectangles,
a stochastic morphism $\kappa\colon S\times X\rightarrow Y$ is a conditional law of $Y$ given $X$ for $\pi$,
that is, satisfies \eqref{eq: CL appendix}, if and only if
\begin{equation}\label{eq: CL integral}
	\int_A \kappa(s,x)(B)\,\pi_X(s)(dx)=\pi(s)(A\times B),\qquad A\in\sB_X,\quad B\in\sB_Y,
\end{equation}
for any $s\in S$.

\begin{proof}[Proof of \cref{thm: main}~\textit{(i)}]
	Let $\kappa$ and $\kappa'$ be conditional laws of $Y$ given $X$ for $\pi$, and fix $s\in S$.
	By \eqref{eq: CL integral}, for any $B\in\sB_Y$ we have
	$\int_A \kappa(s,x)(B)\,\pi_X(s)(dx)=\int_A \kappa'(s,x)(B)\,\pi_X(s)(dx)$ for any $A\in\sB_X$,
	hence $\kappa(s,x)(B)=\kappa'(s,x)(B)$ for $x$ outside a $\pi_X(s)$-null set $N_B$.
	Since $Y$ is countably generated, there exists a countable $\pi$-system $\cG$ generating $\sB_Y$,
	namely the finite intersections of a countable family of generators.
	For $x$ outside the $\pi_X(s)$-null set $\bigcup_{B\in\cG}N_B$,
	the probability measures $\kappa(s,x)$ and $\kappa'(s,x)$ agree on $\cG$,
	hence on $\sB_Y$ by the $\pi$-$\lambda$ theorem.
\end{proof}

The rest of this appendix is devoted to the proof of \cref{thm: main}~\textit{(ii)}.
For $t\in[0,1]$ and a subset $Y\subset[0,1]$, we let $Y^{\leq t}\coloneqq Y\cap[0,t]$.
We first prove the following lemma.

\begin{lemma}\label{lem: construction}
	Let $Y$ be a Borel subset of $[0,1]$, and let $S$ and $X$ be measurable spaces
	such that $X$ is countably generated.  Let
	\[
		\pi\colon S\rightarrow X\times Y
	\]
	be a stochastic morphism, and let $\pi_X\coloneqq\pr_X\circ\pi$ be the marginal.
	Then for any $r\in \bbQ\cap[0,1]$, there exist jointly
	measurable functions $u_r\colon S\times X\rightarrow\R$
	satisfying $u_r(s,x)\leq u_{r'}(s,x)$ if $r\leq r'$, and $u_1(s,x)\equiv 1$,
	such that
	\[
		\int_A u_r(s,x) \pi_X(s)(dx)=\pi(s)(A\times Y^{\leq r})
	\]
	for any $A\in\sB_X$.
\end{lemma}

\begin{proof}
	Since $X$ is countably generated, there exist measurable subsets $E_1$, $E_2$, $\ldots$ of $X$
	such that $\sB_X=\sigma(\{E_i\}_{i=1}^\infty)$.  For integers $n\geq 1$ and $\ve=(\ve_1,\ldots,\ve_n)\in\{0,1\}^n$,
	we let $E^\ve\coloneqq E^{\ve_1}\cap\cdots\cap E^{\ve_n}$, where $E^{\ve_i}=E_i$ if $\ve_i=1$
	and  $E^{\ve_i}=X\setminus E_i$ if $\ve_i=0$.  Then $X=\coprod_{\ve\in\{0,1\}^n}E^\ve$ gives
	a partition of $X$.  Let $\cA^n\coloneqq\{E^\ve\}_{\ve\in\{0,1\}^n}$ be the set of atoms of the partition given by $E_1,\ldots,E_n$,
	and let $\sB^n\coloneqq\sigma(\cA^n)=\sigma(\{E_i\}_{i=1}^n)$.

	For $n\geq 1$ and $r\in\bbQ\cap[0,1]$, set
	\[
		u^n_r(s,x)\coloneqq
			\sum_{E^\ve\in\cA^n} \one_{E^\ve}(x)\frac{\pi(s)(E^\ve\times Y^{\leq r})}{\pi_X(s)(E^\ve)},
	\]
	where $\one_{E^\ve}$ is the indicator function of $E^\ve$
	and
	\begin{equation}
		\frac{\pi(s)(E^\ve\times Y^{\leq r})}{\pi_X(s)(E^\ve)}\coloneqq
		\begin{cases}
			\frac{\pi(s)(E^\ve\times Y^{\leq r})}{\pi_X(s)(E^\ve)}  &\pi_X(s)(E^\ve)>0\\
			r  &  \pi_X(s)(E^\ve)=0.
		\end{cases}
	\end{equation}
	By construction, $u^n_r(s,x)$ is a measurable function on $S\times X$ with values in $[0,1]$
	and is non-decreasing with respect to $r$.  In addition, $u^n_1(s,x)\equiv 1$.

	Fix $s\in S$ and $r\in\bbQ\cap[0,1]$.  We next prove that the functions $\{u^n_r(s,\cdot)\}_{n\geq1}$
	form a martingale on $X$ for the probability measure $\pi_X(s)$ and
	filtration $\{\sB^n\}_{n\geq 1}$.
	For any $E^\ve\in\cA^n$, we have
	\[
		E^\ve = (E^\ve \cap E_{n+1})\amalg (E^\ve \cap E^c_{n+1}),\qquad
		E^c_{n+1}\coloneqq X\setminus E_{n+1}.
	\]
	By definition of $u^{n+1}_r(s,x)$, we have
	\begin{align*}
		\int_{E^\ve} u^{n+1}_r(s,x)\pi_X(s)(dx)&=
		\int_{E^\ve\cap E_{n+1}} u^{n+1}_r(s,x)\pi_X(s)(dx)+\int_{E^\ve\cap E^c_{n+1}} u^{n+1}_r(s,x)\pi_X(s)(dx)\\
		& =\frac{\pi(s)((E^\ve\cap E_{n+1})\times Y^{\leq r})}{\pi_X(s)(E^\ve\cap E_{n+1})} \pi_X(s)(E^\ve\cap E_{n+1}) \\
		 &\qquad+\frac{\pi(s)((E^\ve\cap E^c_{n+1})\times Y^{\leq r})}{\pi_X(s)(E^\ve\cap E^c_{n+1})} \pi_X(s)(E^\ve\cap E^c_{n+1}).
	\end{align*}
	If $\pi_X(s)(E^\ve)>0$, then the above sum equals $\pi(s)(E^\ve\times Y^{\leq r})$,
	and if $\pi_X(s)(E^\ve)=0$, then the above sum equals $0$.
	Thus we have
	\begin{align*}
		\int_{E^\ve} u^{n+1}_r(s,x)\pi_X(s)(dx)&=\int_{E^\ve} u^{n}_r(s,x)\pi_X(s)(dx)
	\end{align*}
	for any $E^\ve\in\cA^n$, where the integral on both sides is zero
	for $E^\ve$ such that $\pi_X(s)(E^\ve)=0$.
	By definition of the conditional expectation, this shows that
	\[
		\bbE_{\pi_X(s)}[u^{n+1}_r(s,\cdot)\mid\sB^n]=u^n_r(s,\cdot),
	\]
	hence $\{u^n_r(s,\cdot)\}_{n\geq1}$
	form a martingale on $X$ for $\pi_X(s)$ and filtration $\{\sB^n\}_{n\geq 1}$.

	By the bounded martingale convergence theorem,
	for each $s\in S$, the sequence
	$\{u^n_r(s,\cdot)\}_{n\geq1}$ converges $\pi_X(s)$-almost everywhere
	and in $L^1(X,\pi_X(s))$ to a function $u_r(s,\cdot)$.  For any $(s,x)\in S\times X$, let
	\[
		u_r(s,x)\coloneqq\limsup_{n\rightarrow\infty}u_r^n(s,x).
	\]
	Then $u_r(s,x)$ is a jointly measurable function on $S\times X$
	which coincides with the pointwise and $L^1$-limits $\pi_X(s)$-almost everywhere.
	We have
	\[
		\Bigl|\int_A u^n_r(s,x)\pi_X(s)(dx)-\int_A u_r(s,x)\pi_X(s)(dx)\Bigr|
		\leq\int_X\abs{u^n_r(s,x)- u_r(s,x)}\pi_X(s)(dx),
	\]
	thus
	\[
		\lim_{n\rightarrow\infty}\int_A u^n_r(s,x)\pi_X(s)(dx)=\int_A u_r(s,x)\pi_X(s)(dx).
	\]
	For $A\in\sB^n$,
	by definition of $u^m_r(s,x)$, for $m\geq n$, we have
	$
		\int_{A} u^{m}_r(s,x)\pi_X(s)(dx)=\pi(s)(A\times Y^{\leq r}).
	$
	This shows that $u_r(s,x)$ satisfies
	\[
		\int_A u_r(s,x)\pi_X(s)(dx)=\pi(s)(A\times Y^{\leq r})
	\]
	for any $A\in\bigcup_{n\geq 1}\sB^n$.
	The measures are probability measures and hence in particular $\sigma$-finite.
	The equality extends to any $A\in\sB_X=\sigma(\bigcup_{n\geq 1}\sB^n)$,
	since  any finite measures that agree on an algebra necessarily agree on the generated
	$\sigma$-algebra.
	This completes the proof of our assertion.
\end{proof}

\begin{proof}[Proof of \cref{thm: main}~\textit{(ii)}]
	Since we have assumed that $Y$ is a standard Borel space, it is isomorphic as a measurable space to a Borel subset of $[0,1]$.
	We identify $Y$ with its image and assume that $Y\subset[0,1]$.
	By \cref{lem: construction}, for any $r\in\bbQ\cap[0,1]$, there exist
	jointly measurable functions $u_r\colon S\times X\rightarrow\R$
	satisfying $u_r(s,x)\leq u_{r'}(s,x)$ if $r\leq r'$ and $u_1(s,x)\equiv 1$,
	such that
	\[
		\int_A u_r(s,x) \pi_X(s)(dx)=\pi(s)(A\times Y^{\leq r})
	\]
	for any $A\in\sB_X$.
	Let
	\[
		F_t(s,x)\coloneqq\inf_{r\in\bbQ\cap(t,1]} u_r(s,x),\qquad t\in[0,1)
	\]
	and $F_1(s,x)\equiv 1$.  For fixed $(s,x)\in S\times X$, the function $t\mapsto F_t(s,x)$ on
	$[0,1]$ is a non-decreasing right continuous function with $F_1(s,x)\equiv1$.
	By the theorem of Lebesgue--Stieltjes, $t\mapsto F_t(s,x)$ is a distribution function of a unique
	probability measure $\kappa(s,x)$ on $[0,1]$, such that
	\[
		\kappa(s,x)([0,t])=F_t(s,x)
	\]
	for any $t\in[0,1]$.  For a fixed $t$, the function $F_t(s,x)$ is a countable infimum
	of jointly measurable functions, hence is jointly measurable.
	For $t\in[0,1)$, take a sequence of rational points $r_m\downarrow t$.
	Then $u_{r_m}(s,x)\rightarrow F_t(s,x)$ by monotonicity, so by the dominated convergence theorem
	and continuity from above of $\pi(s)$, we see that
	\[
		\int_A \kappa(s,x)([0,t])\pi_X(s)(dx)=\pi(s)(A\times Y^{\leq t}),\qquad t\in[0,1], \quad A\in\sB_X.
	\]
	The system
	\[
		\sD=\{\tilde B\in \sB_{[0,1]}\mid \text{$\kappa(s,x)(\tilde B)$ jointly measurable in $(s,x)$}\}
	\]
	is a $\lambda$-system containing  the $\pi$-system $\{[0,t]\mid t\in[0,1]\}$ generating $\sB_{[0,1]}$.
	  Hence again by the  $\pi$-$\lambda$ theorem, we have $\sD=\sB_{[0,1]}$.
	This shows that $\kappa(s,x)(\tilde B)$ is jointly measurable in $(s,x)$ for any $\tilde B\in\sB_{[0,1]}$.
	We next extend the above equality from $[0,t]$ to any measurable $\tilde B\subset[0,1]$.
	Fix $s\in S$ and $A\in\sB_X$, and consider the finite measures on $[0,1]$ given by
	\begin{align*}
		\int_A&\kappa(s,x)(\tilde B)\pi_X(s)(dx),&
		\pi(s)&(A\times(Y\cap \tilde B))
	\end{align*}
	for any $\tilde B\in\sB_{[0,1]}$.
	They agree on the $\pi$-system $\{[0,t]\mid t\in[0,1]\}$ generating $\sB_{[0,1]}$,
	hence by the $\pi$-$\lambda$ theorem coincide as finite measures on $[0,1]$.
	Thus we have
	\[
		\int_A \kappa(s,x)(\tilde B)\pi_X(s)(dx)=\pi(s)(A\times(Y\cap\tilde B)),\qquad A\in\sB_X
	\]
	for any $\tilde B\in\sB_{[0,1]}$.
	In particular,
	\[
		\int_X \kappa(s,x)(Y)\pi_X(s)(dx)=\pi(s)(X\times Y)=1.
	\]
	Since $\kappa(s,x)(Y)\leq 1$, the set $N_s\coloneqq\{x \in X \mid \kappa(s,x)(Y)<1\}$ is
	$\pi_X(s)$-null.
	Moreover, since $Y\in\sD$, the function $\kappa(s,x)(Y)$ is jointly measurable on $S\times X$,
	so the set $N\coloneqq\{(s,x) \in S\times X \mid \kappa(s,x)(Y)<1\}$ is measurable in $S\times X$.
	Fix $y_0\in Y$ and modify $\kappa$ so that $\kappa(s,x)=\delta_{y_0}$ on $N$.
	This modification
	does not change the fact that $\kappa(s,x)(\tilde B)$ is jointly measurable on $S\times X$, and
	ensures that the measure $\kappa(s,x)$ is supported on $Y$ for any $(s,x)\in S\times X$.
	Since $\sB_Y\subset\sB_{[0,1]}$, for any $B\in\sB_Y$,
	the function $\kappa(s,x)(B)$ is jointly measurable on $S\times X$,
	hence $\kappa$ gives a stochastic morphism $\kappa\colon S\times X\rightarrow Y$.
	Moreover,
	\[
		\int_A \kappa(s,x)(B)\pi_X(s)(dx)=\pi(s)(A\times B),\qquad A\in\sB_X, \quad B\in\sB_Y,
	\]
	which is \eqref{eq: CL integral}.  Hence $\kappa$ is a conditional law of $Y$ given $X$ for $\pi$.
	The uniqueness follows from \textit{(i)}, since a standard Borel space is countably generated.
	This completes the proof of \cref{thm: main}~\textit{(ii)}.
\end{proof}

%
%
%
\section{Expectation and Conditional Expectation}\label{sec: appendix B}
%
%
%

In this appendix, we define the expectations and conditional expectations
as pull-backs of functions with respect to conditional laws.
Throughout this appendix, we assume that the parameter space $S=\bsone\coloneqq\{*\}$,
and will identify $\mu\in\sP(X)$ with a stochastic morphism $\mu\colon\bsone\rightarrow X$
mapping $*$ to $\mu$.  Here, $\sP(X)$ is the space of probability measures on $X$.
For a stochastic morphism $u\colon X\rightarrow Y$, we write $u(\mu)\coloneqq u\circ\mu$
for the probability measure on $Y$ given by $u(\mu)(B)=\int_Xu(x)(B)\mu(dx)$ for $B\in\sB_Y$.

For any measurable space $X$, denote by $\cB(X)$ the space of bounded measurable functions on $X$.
The space $\cB(X)$ is a Banach space with respect to the supremum
norm $\dabs{f}_\infty\coloneqq\sup_{x\in X}\abs{f(x)}$
for any $f\in\cB(X)$.  We define the \emph{expectation} of $f$ with respect to $\mu\in\sP(X)$ by
\[
	\bbE_\mu[f]\coloneqq\int_Xf(x)\mu(dx).
\]

\begin{definition}\label{def: pullback}
	Let $u\colon X\rightarrow Y$ be a stochastic morphism.
	For any $g\in\cB(Y)$, we define the \emph{pull-back} $u^*g\in\cB(X)$
	to be the bounded measurable function
	\[
		u^*g(x)\coloneqq\int_Yg(y) u(x)(dy),\qquad x\in X.
	\]
\end{definition}

Note that if $u\colon X\rightarrow Y$ is a deterministic morphism, then
\[
	u^*g(x)=\int_Y g(y)\delta_{u(x)}(dy)=g(u(x))=(g\circ u)(x)
\]
for any $g\in\cB(Y)$ and $x\in X$, hence we have $u^* g=g\circ u$.

\begin{proposition}\label{prop: contravariant}
	Let $u\colon X\rightarrow Y$ be a stochastic morphism.
	Then the pull-back
	\[
		u^*\colon\cB(Y)\rightarrow\cB(X)
	\]
	is a continuous map of operator norm $\leq 1$.
	In addition, for stochastic morphisms $u\colon X\rightarrow Y$
	and $v\colon Y\rightarrow Z$,
	the pull-back of the composition satisfies
	\begin{equation}\label{eq: composition}
		(v\circ u)^*= u^*\circ v^*\colon\cB(Z)\rightarrow\cB(X).
	\end{equation}
\end{proposition}

\begin{proof}
	Given a stochastic morphism $u\colon X\rightarrow Y$
	and $g\in\cB(Y)$, the pull-back $u^*g$
	satisfies
	\[
		\abs{u^*g(x)}\leq\int_Y\abs{g(y)} u(x)(dy)\leq\dabs{g}_\infty,\qquad x\in X.
	\]
	This shows that $u^*\colon\cB(Y)\rightarrow\cB(X)$ has operator norm $\leq1$,
	hence gives a continuous map of Banach spaces.
	In addition, if $u\colon X\rightarrow Y$ and $v\colon Y\rightarrow Z$ are stochastic morphisms,
	then
	\begin{align*}
		u^*(v^*h)(x)&=\int_Y v^*h(y) u(x)(dy)=\int_Z\int_Y h(z) v(y)(dz)u(x)(dy)\\
		&=\int_Z h(z) (v\circ u)(x)(dz)=(v\circ u)^* h(x),
	\end{align*}
	which proves \eqref{eq: composition}.
\end{proof}

We next extend the pull-back to $L^p$-spaces.
For any probability measure $\mu$ on $X$ and $p\geq1$, let $L^p(X,\mu)$ be the $L^p$ space consisting of
equivalence classes of $p$-integrable functions on $X$ with respect to the probability measure $\mu$,
where two functions are equivalent if they agree $\mu$-almost everywhere.
The $L^p$-norm is defined by
\[
	\dabs{f}_{L^p(X,\mu)}\coloneqq\biggl(\int_X\abs{f(x)}^p\mu(dx)\biggr)^{1/p}.
\]
The space of bounded measurable functions $\cB(X)$
is dense in $L^p(X,\mu)$ for the $L^p$-norm.

\begin{lemma}\label{lem: Lp}
	Let $\mu$ be a probability measure on $X$, and let $u\colon X\rightarrow Y$ be a stochastic morphism.
	Then for any $p\geq1$, the pull-back extends uniquely to a continuous map
	\[
		u^*\colon L^p(Y,u(\mu))\rightarrow L^p(X,\mu)
	\]
	of operator norm $\leq1$.
\end{lemma}

\begin{proof}
	By definition of $u(\mu)$, for any $h\in\cB(Y)$ we have
	\[
		\bbE_\mu[u^*h]=\int_X\int_Yh(y)u(x)(dy)\mu(dx)=\int_Yh(y)u(\mu)(dy)=\bbE_{u(\mu)}[h].
	\]
	Let $g\in\cB(Y)$.  Since $u(x)$ is a probability measure on $Y$ for any $x\in X$,
	Jensen's inequality for the convex function $t\mapsto\abs{t}^p$ gives
	\[
		\abs{u^*g(x)}^p=\Bigl|\int_Yg(y)u(x)(dy)\Bigr|^p\leq\int_Y\abs{g(y)}^pu(x)(dy)=u^*\abs{g}^p(x),
		\qquad x\in X,
	\]
	hence applying the above equality to $h=\abs{g}^p$, we have
	\[
		\dabs{u^*g}^p_{L^p(X,\mu)}=\bbE_\mu[\abs{u^*g}^p]\leq\bbE_\mu[u^*\abs{g}^p]
		=\bbE_{u(\mu)}[\abs{g}^p]=\dabs{g}^p_{L^p(Y,u(\mu))}.
	\]
	In particular, if $g=0$ holds $u(\mu)$-almost everywhere, then $u^*g=0$ holds $\mu$-almost everywhere,
	hence $u^*$ induces a map of operator norm $\leq1$ from the image of $\cB(Y)$ in $L^p(Y,u(\mu))$
	to $L^p(X,\mu)$.  Since $\cB(Y)$ is dense in $L^p(Y,u(\mu))$, this map extends uniquely to a continuous map
	$u^*\colon L^p(Y,u(\mu))\rightarrow L^p(X,\mu)$ of operator norm $\leq1$.
\end{proof}

The value $u^*g(x)=\bbE_{u(x)}[g]$ is the expectation of $g$ for the probability measure $u(x)$ on $Y$.
Hence we will often write
\[
	\bbE_u[g]\coloneqq u^*g
\]
and refer to the pull-back $u^*g$ as the expectation of $g$ with respect to $u$.
Consider a probability measure $\pi\in\sP(X\times Y)$, which we also
view as a stochastic morphism $\pi\colon\bsone\rightarrow X\times Y$.
In what follows, we assume that
the conditional laws $\pi(X\mid Y)\colon Y\rightarrow X$ and $\pi(Y\mid X)\colon X\rightarrow Y$ exist.
The expectations of conditional laws satisfy the following adjoint property.

\begin{lemma}\label{lem: CE}
	Let $\pi\in\sP(X\times Y)$.
	Then for any $f\in\cB(X)$ and $g\in\cB(Y)$, we have
	\begin{equation}\label{eq: adjoint}
		\bbE_{\pi_X}[f \bbE_{\pi(Y\mid X)}[g]] =
		\bbE_\pi[fg]=\bbE_{\pi_Y}[\bbE_{\pi(X\mid Y)}[f] g  ],
	\end{equation}
	where $\pi_X$ and $\pi_Y$ are the $X$- and $Y$-marginals of $\pi$,
	and $fg\in\cB(X\times Y)$ denotes the function $(x,y)\mapsto f(x)g(y)$
	on $X\times Y$.
\end{lemma}

\begin{proof}
	For $S=\bsone$, the conditional law is characterized by
	$\pi=(\pr_X,\pi(Y\mid X))\circ\pi_X$.
	Hence
	\begin{align*}
		\bbE_\pi[fg]&=\int_{X\times Y}f(x)g(y)\int_X\delta_{x'}(dx)\pi(Y\mid X)(x)(dy)\,\pi_X(dx')\\
		&=\int_{X}f(x)\int_Y g(y)\pi(Y\mid X)(x)(dy)\pi_X(dx)=\bbE_{\pi_X}[f\bbE_{\pi(Y\mid X)}[g]],
	\end{align*}
	which proves the first equality.
	The second equality follows similarly by exchanging the roles of $X$ and $Y$.
\end{proof}

The spaces $L^2(X,\pi_X)$ and $L^2(Y,\pi_Y)$ are Hilbert spaces
with inner products
\begin{align*}
	\pair{f,f'}_{\pi_X}&\coloneqq\bbE_{\pi_X}[ff'],&
	\pair{g,g'}_{\pi_Y}&\coloneqq\bbE_{\pi_Y}[gg']
\end{align*}
for $f,f'\in L^2(X,\pi_X)$ and $g,g'\in L^2(Y,\pi_Y)$.  Then \eqref{eq: adjoint} extended to $L^2$ by continuity gives
\[
		\pair{f , \bbE_{\pi(Y\mid X)}[g]}_{\pi_X} =
		\pair{\bbE_{\pi(X\mid Y)}[f] ,g}_{\pi_Y},
\]
showing that $\bbE_{\pi(Y\mid X)}$ and $\bbE_{\pi(X\mid Y)}$ are adjoint to each other as
operators between Hilbert spaces.

For $\pi\in\sP(X\times Y)$, the \emph{conditional expectation}  $\bbE_\pi[f\mid Y]$ of $f\in L^1(X,\pi_X)$
given $Y$ for the joint law $\pi$ is characterized as the unique function in $L^1(Y,\pi_Y)$ satisfying
\begin{equation}\label{eq: characterization}
	\bbE_{\pi_Y}[\bbE_\pi[f\mid Y]g]=\bbE_\pi[fg],\qquad g\in\cB(Y).
\end{equation}
\cref{lem: CE} gives the following.

\begin{lemma}\label{lem: CE CL}
	Let $\pi\in\sP(X\times Y)$.
	Then for any $f\in L^1(X,\pi_X)$, we have
	\[
		\bbE_\pi[f\mid Y]=\bbE_{\pi(X\mid Y)}[f]\in L^1(Y,\pi_Y).
	\]
\end{lemma}

\begin{proof}
	By \cref{lem: CE}, $\bbE_{\pi(X\mid Y)}[f]$ satisfies the characterization \eqref{eq: characterization} of the
	conditional expectation for $f\in\cB(X)$, and for $f\in L^1(X,\pi_X)$ by continuity,
	since both sides of \eqref{eq: characterization} are continuous in $f$ for the $L^1$-norm
	by $\abs{\bbE_\pi[fg]}\leq\dabs{g}_\infty\dabs{f}_{L^1(X,\pi_X)}$ and \cref{lem: Lp}.
	For the uniqueness, suppose $h_1,h_2\in L^1(Y,\pi_Y)$ both satisfy $\bbE_{\pi_Y}[h_i\,g]=\bbE_\pi[fg]$ for all $g\in\cB(Y)$. Taking $g=\one_B$ for $B\in\sB_Y$ gives
\[
	\int_Bh_1\,d\pi_Y=\bbE_\pi[f\one_B]=\int_Bh_2\,d\pi_Y\qquad\text{for all }B\in\sB_Y,
\]
so $\int_B(h_1-h_2)\,d\pi_Y=0$ for every $B\in\sB_Y$. Choosing $B=\{h_1>h_2\}$, the integrand is strictly positive on $B$, forcing $\pi_Y(B)=0$.  Similarly for $\{h_1<h_2\}$. Hence $h_1=h_2$ holds $\pi_Y$-almost everywhere, i.e.\ they are the same element of $L^1(Y,\pi_Y)$.
\end{proof}

The tower property of the conditional law (\cref{thm: tower}), in the form of the
chain rule \cref{lem: properties}~\ref{CL: C}, gives the tower property of conditional expectations.
Consider $\pi\in\sP(X\times Y\times Z)$.  The conditional laws and expectations below
are taken with respect to the marginals of $\pi$.

\begin{proposition}\label{prop: tower}
	Let $\pi\in\sP(X\times Y\times Z)$, and suppose the conditional laws
	$\pi(Y\mid X)$ and $\pi(Z\mid X\times Y)$ exist.
	Then $\pi(Z\mid X)$ exists, and for any $f\in\cB(Z)$, we have
	\[
		\bbE_\pi[\bbE_\pi[f\mid X\times Y]\bigm| X]=\bbE_\pi[f\mid X].
	\]
\end{proposition}

\begin{proof}
	By the chain rule \cref{lem: properties}~\ref{CL: C} for $S=\bsone$,
	the conditional law $\pi(Z\mid X)$ exists and satisfies
	$\pi(Z\mid X)=\pi(Z\mid X\times Y)\circ\pi(Y\mid X)^\wedge$.
	By \cref{lem: CE CL}, we have
	$\bbE_\pi[f\mid X\times Y]=\bbE_{\pi(Z\mid X\times Y)}[f]$ and $\bbE_\pi[f\mid X]=\bbE_{\pi(Z\mid X)}[f]$.
	For $h\in\cB(X\times Y)$, the characterization \eqref{eq: characterization} of $\bbE_\pi[h\mid X]$
	and $\pi_{X\times Y}=\pi(Y\mid X)^\wedge\circ\pi_X$ of \eqref{eq: CL aug} give
	\[
		\bbE_{\pi_X}[\bbE_\pi[h\mid X]\,g]=\bbE_{\pi_{X\times Y}}[h\,(g\circ\pr_X)]
		=\bbE_{\pi_X}[\bbE_{\pi(Y\mid X)^\wedge}[h]\,g],\qquad g\in\cB(X),
	\]
	hence $\bbE_\pi[h\mid X]=\bbE_{\pi(Y\mid X)^\wedge}[h]$.
	Hence by \eqref{eq: composition}, we have
	\begin{align*}
		\bbE_\pi[f\mid X]&=\pi(Z\mid X)^*f
		=(\pi(Y\mid X)^\wedge)^*(\pi(Z\mid X\times Y)^*f)
		=\bbE_\pi[\bbE_\pi[f\mid X\times Y]\bigm| X]
	\end{align*}
	as desired.
\end{proof}

Consider random variables $\bX$, $\bY$ on a probability space
$(\Omega,\sF,\bbP)$ with values in measurable spaces $(X,\sB_X)$, $(Y,\sB_Y)$.
Let $\pi$ be the joint law of $\bX$ and $\bY$ on the canonical probability space 
$(X\times Y,\sB_X\otimes\sB_Y)$.
The coordinate projections $\pr_X$ and $\pr_Y$ on the probability space $(X\times Y,\pi)$
have the same joint law as $\bX$ and $\bY$, so that
\[
	\bbE_\bbP[f(\bX)]=\bbE_\pi[f\circ\pr_X]=\bbE_{\pi_X}[f]
\]
for any bounded measurable $f\colon X\rightarrow\R$.
Thus, our expectation notation $\bbE_\pi$ agrees with classical expectation under this canonical
realization.

%
%
%
\section{A Dictionary of Notations}\label{sec: appendix C}
%
%
%

In this appendix, we compare the notation of this article with the classical notation of probability theory
and with the notation of Markov categories \cite{F20}, in which morphisms are drawn as string diagrams
read from bottom to top, with $\mathrm{copy}_X\colon X\rightarrow X\otimes X$ and $\mathrm{del}_X\colon X\rightarrow I$
the copy and delete morphisms, and $\mathrm{copy}^{(3)}_X\colon X\rightarrow X\otimes X\otimes X$ the threefold copy.
In the classical column, $X$, $Y$ and $Z$ denote random variables on a probability space $(\Omega,\sF,P)$.
In the other two columns, they denote measurable spaces or objects.
Classical conditional probabilities and conditionals in Markov categories are unique only almost surely,
as is the conditional law of this article.

\begingroup
\renewcommand{\arraystretch}{1.6}
\small
\begin{longtable}{|>{\raggedright\arraybackslash}p{0.15\textwidth}|>{\raggedright\arraybackslash}p{0.25\textwidth}|>{\raggedright\arraybackslash}p{0.25\textwidth}|>{\raggedright\arraybackslash}p{0.25\textwidth}|}
\hline
 & \textbf{Classical} & \textbf{Markov category} & \textbf{This article}\\
\hline
\endhead
Law of $X$
& $P_X=P\circ X^{-1}$
& state $p\colon I\rightarrow X$\quad
\begin{tikzpicture}[baseline={(0,0.4)},scale=0.8]
	\draw (0,0.5)--(0,1.1) node[above]{\scriptsize $X$};
	\draw (-0.35,0) rectangle (0.35,0.5); \node at (0,0.25){\scriptsize $p$};
\end{tikzpicture}
& $\mu\colon\bsone\rightarrow X$\\
\hline
Family of laws
& $(P_\theta)_{\theta\in\Theta}$
& $p\colon A\rightarrow X$, with $A$ an input wire
& $\mu\colon S\rightarrow X$, with $S$ the parameter space\\
\hline
Markov kernel
& $k(x,dy)$
& $f\colon X\rightarrow Y$
& $h\colon X\rightarrow Y$\\
\hline
Law of $f(X)$
& $P\circ (f\circ X)^{-1}$
& $f\circ p$
& $f\circ\mu$\\
\hline
Joint law of $X$ and $Y$
& law of $(X,Y)$
& $\psi\colon I\rightarrow X\otimes Y$\quad
\begin{tikzpicture}[baseline={(0,0.4)},scale=0.8]
	\draw (-0.25,0.5)--(-0.25,1.1) node[above]{\scriptsize $X$};
	\draw (0.25,0.5)--(0.25,1.1) node[above]{\scriptsize $Y$};
	\draw (-0.5,0) rectangle (0.5,0.5); \node at (0,0.25){\scriptsize $\psi$};
\end{tikzpicture}
& $\pi\colon S\rightarrow X\times Y$\\
\hline
Product with parameter
& $P_\theta\otimes Q_\theta$ for each $\theta$
& $(p\otimes q)\circ\mathrm{copy}_A$ for $p\colon A\rightarrow X$, $q\colon A\rightarrow Y$
& $(p,q)\colon S\rightarrow X\times Y$\\
\hline
Parallel product 
& $k(x,\cdot)\otimes l(y,\cdot)$
& $f\otimes g$
& $f\times g$\\
\hline
Composition along a shared parameter
& $\int_V r(s,v,\cdot)\,q(s,u,dv)$
& $r\circ(\id_S\otimes q)\circ(\mathrm{copy}_S\otimes\id_U)$
& $r\circ q^\wedge$\\
\hline
Marginal
& $P_X$
& $(\id_X\otimes\mathrm{del}_Y)\circ\psi$\quad
\begin{tikzpicture}[baseline={(0,0.4)},scale=0.8]
	\draw (-0.25,0.5)--(-0.25,1.1) node[above]{\scriptsize $X$};
	\draw (0.25,0.5)--(0.25,0.85);
	\draw (0.1,0.85)--(0.4,0.85); \draw (0.15,0.93)--(0.35,0.93); \draw (0.2,1.01)--(0.3,1.01);
	\draw (-0.5,0) rectangle (0.5,0.5); \node at (0,0.25){\scriptsize $\psi$};
\end{tikzpicture}
& $\pi_X=\pr_X\circ\pi$\\
\hline
Conditional law of $Y$ given $X$
& $P(Y\in\cdot\mid X=x)$
& conditional $\psi_{|X}\colon X\otimes A\rightarrow Y$ of $\psi\colon A\rightarrow X\otimes Y$
& $\pi(Y\mid X)\colon S\times X\rightarrow Y$\\
\hline
Conditional law with the condition copied
& $\delta_x\otimes P(Y\in\cdot\mid X=x)$
& $(\id_X\otimes\psi_{|X})\circ(\mathrm{copy}_X\otimes\id_A)\colon X\otimes A\rightarrow X\otimes Y$
& $\pi(X\times Y\mid X)=(\pr_X,\pi(Y\mid X))$, with $\pi=\pi(X\times Y\mid X)\circ(\pi_X)^\wedge$, see \eqref{eq: CL copy}\\
\hline
Disintegration
& $P_{(X,Y)}(dx\,dy)=P_X(dx)$\newline \hspace*{1em}\mbox{$P(Y\in dy\mid X=x)$}
&
\begin{tikzpicture}[baseline={(0,0.9)},scale=0.8]
	\draw (-0.4,0) rectangle (0.4,0.5); \node at (0,0.25){\scriptsize $\psi_X$};
	\draw (0,0.5)--(0,0.8); \fill (0,0.8) circle (0.06);
	\draw (0,0.8)--(-0.45,1.15)--(-0.45,2.2) node[above]{\scriptsize $X$};
	\draw (0,0.8)--(0.45,1.15)--(0.45,1.3);
	\draw (0.05,1.3) rectangle (0.85,1.8); \node at (0.45,1.55){\scriptsize $\psi_{|X}$};
	\draw (0.45,1.8)--(0.45,2.2) node[above]{\scriptsize $Y$};
\end{tikzpicture}
\;$=$\;
\begin{tikzpicture}[baseline={(0,0.9)},scale=0.8]
	\draw (-0.25,0.5)--(-0.25,2.2) node[above]{\scriptsize $X$};
	\draw (0.25,0.5)--(0.25,2.2) node[above]{\scriptsize $Y$};
	\draw (-0.5,0) rectangle (0.5,0.5); \node at (0,0.25){\scriptsize $\psi$};
\end{tikzpicture}
& $\pi^\wedge=\pi(Y\mid X)^\wedge\circ(\pi_X)^\wedge$, that is,
\newline
\hspace*{0.5em}$\xymatrix@R=1.4em@C=1.6em{
	S \ar[dd]_{\pi^\wedge} \ar[dr]^{(\pi_X)^\wedge} & \\
	& S\times X \ar[dl]^(.4){\;\pi(Y\mid X)^\wedge}\\
	S\times X\times Y &
}$
\newline commutes, see \eqref{eq: CL aug}\\
\hline
Tower property for laws
& no expression without auxiliary kernels; for densities, $p(z\mid x,y)=\frac{p(y,z\mid x)}{p(y\mid x)}$
& every $(\psi_{|X})_{|Y}$ is a conditional $\psi_{|X\otimes Y}$, up to the order of inputs \cite[Lemma~11.11]{F20}
& $\pi(Y\times Z\mid X)(Z\mid Y)=$\newline $\pi(Z\mid X\times Y)$ (\cref{thm: tower})\\
\hline
Conditional expectation
& $\bbE[f(Y)\mid X]$
& ---
& $\bbE_\pi[f\mid X]=\pi(Y\mid X)^*f$ (Appendix \ref{sec: appendix B})\\
\hline
Conditional independence
& $X\perp Y\mid W$
& for $\psi\colon I\rightarrow W\otimes X\otimes Y$, $\psi=(\id_W\otimes f\otimes g)\circ\mathrm{copy}^{(3)}_W\circ\varphi$ for some $f\colon W\rightarrow X$, $g\colon W\rightarrow Y$, $\varphi\colon I\rightarrow W$ \cite[Definition~12.1]{F20}
& $X\perp_\pi Y\mid W$, i.e.\ $\pi(X\times Y\mid W)=$\newline $(\pi(X\mid W),\pi(Y\mid W))$ (\cref{def: CI}), equivalently $\pi(Y\mid W\times X)=$\newline $\pi(Y\mid W)$ (\cref{lem: CI})\\
\hline
Bayes' rule
& $p(x\mid y)\propto p(y\mid x)\,p(x)$
& Bayesian inverse $f^\dagger_p\colon Y\rightarrow X$
& $h^\mu\colon S\times Y\rightarrow X$, with parameterized prior $\mu\colon S\rightarrow X$\\
\hline
Bayesian filter
& $p(x_k\mid y_1,\ldots,y_k)$
& Bayes filter of \cite{F25}, without feedback
& $\bspi(X_k\mid Y^{[k]})=h_k^{\mu_k}$, with feedback (\cref{thm: BF})\\
\hline
Affine Gaussian kernel
& $Y=AX+b+\varepsilon$, $\varepsilon\sim\sN(0,R)$ independent of $X$
& morphism of $\mathbf{Gauss}$, with $A$, $b$, $R$ constant
& $\sN(AX+\bsb,R)$, with $A$, $\bsb$, $R$ depending on $s\in S$\\
\hline
\end{longtable}
\endgroup

\subsection*{Acknowledgment}

The first author would like to thank Emtiyaz Khan for introducing to him \cite{J03},
and for inspiration and passionate discussion through the Bayes-Duality Project,
which helped him to understand the Bayesian way of thinking.
The authors thank Thomas M\"ollenhoff and Yoshino Taro 
for helpful insights.  We  also thank
the members of the Bayes Study Group at Keio University/RIKEN AIP
for discussion, especially Kei Hagihara and Kanji Inui.

\end{document}